\documentclass[final,5p,times,twocolumn,authoryear]{elsarticle}

\usepackage{amsmath,amssymb,amsfonts,mathtools}
\usepackage{amsthm}
\usepackage{algorithm}
\usepackage{algpseudocode}
\usepackage{bm}
\usepackage{booktabs}
\usepackage{graphicx}
\usepackage[hidelinks]{hyperref}
\usepackage{microtype}

\mathtoolsset{showonlyrefs}

\newcommand{\R}{\mathbb{R}}
\newcommand{\cH}{\mathcal{H}}
\newcommand{\diag}{\operatorname{diag}}
\newcommand{\tr}{\operatorname{tr}}
\newcommand{\dd}{\,\mathrm{d}}
\providecommand{\tcite}{\citep}
\providecommand{\doi}[1]{%
  \href{https://doi.org/#1}{doi:\nolinkurl{#1}}%
}

\newtheorem{theorem}{Theorem}[section]
\newtheorem{proposition}[theorem]{Proposition}
\newtheorem{lemma}[theorem]{Lemma}

\theoremstyle{definition}
\newtheorem{assumption}[theorem]{Assumption}

\theoremstyle{remark}
\newtheorem{remark}[theorem]{Remark}

\renewcommand{\tcite}{\citet}

\graphicspath{{figures/}{./}}

\newenvironment{widefigure}[1][htbp]{\begin{figure*}[#1]}{\end{figure*}}

\newcommand{\columnfigurewidth}{\linewidth}

\begin{document}

\begin{frontmatter}

  \title{Controllability Allocation Scores for Linear Infinite-Dimensional Systems}

  \author[imperial]{Daisuke Inoue\corref{cor1}}
  \ead{dinoue@ic.ac.uk}
  \author[tokyo]{Kazuhiro Sato}
  \ead{kazuhiro@mist.i.u-tokyo.ac.jp}

  \cortext[cor1]{Corresponding author.}
  \affiliation[imperial]{%
    organization={Department of Mathematics, Imperial College London},
    addressline={Huxley Building, South Kensington Campus},
    city={London},
    postcode={SW7 2AZ},
    country={United Kingdom}%
  }
  \affiliation[tokyo]{%
    organization={Department of Mathematical Informatics, Graduate School of Information Science and Technology, The University of Tokyo},
    addressline={7-3-1 Hongo, Bunkyo-ku},
    city={Tokyo},
    postcode={113-8656},
    country={Japan}%
  }

  \begin{abstract}
    We define finite-time controllability allocation scores (CASs) as optimal solutions to two convex optimization problems that allocate a fixed actuation budget among prescribed input functions of linear infinite-dimensional continuum systems with selected target outputs.
The objectives are expressed using controllability Gramians for the selected outputs and correspond to maximizing the volume of output changes reachable with unit control energy or minimizing the mean energy required for unit output changes.
To compute these CASs, we discretize the system in space, construct the resulting Gramians, and solve the resulting allocation problems by a domain-aware projected gradient descent method.
We study the behavior of CASs when the actuator-to-target responses of the spatially discretized systems converge to those of the continuum system as the mesh is refined.
If the continuum system is target controllable, we show that the spatially discretized systems are also target controllable on sufficiently fine meshes and that the optimal values of the approximate problems converge to the optimal value for the continuum system.
If the allocation problem for the continuum system has a unique CAS, we show that the CASs obtained from the spatially discretized problems converge to it.
We also derive quantitative bounds on the errors between these CASs when the Gramian map for the continuum system distinguishes allocations.
A two-dimensional advection--diffusion example illustrates how transport dynamics affect target controllability and CASs.
Spatial refinement reveals apparent controllability caused by numerical diffusion and supports the theoretical convergence results.

  \end{abstract}


\end{frontmatter}

\section{Introduction}
\label{sec:introduction}

Controlling spatially distributed systems often requires allocating limited
actuation resources among candidate inputs.  Their states are spatial
fields, such as temperature, concentration, or displacement, that evolve
over time and are represented in infinite-dimensional function spaces
\citep{CurtainZwart1995,TucsnakWeiss2009}.  Examples include room temperature
control with multiple heaters and ventilation outlets
\citep{HeGonzalez2016} and suppression of plate vibrations using
piezoelectric patches \citep{Morris2020}.  Determining an optimal allocation
therefore requires comparing control performance under different allocations
in these infinite-dimensional systems.

To evaluate these allocations, one can take either a model-free or a
model-based approach.  The former uses measured system responses directly,
avoiding the need to identify an explicit model of the dynamics
\citep{BaggioKatewaPasqualetti2019,BaggioBassettPasqualetti2021,vanWaardeEisingTrentelmanCamlibel2020}.
The latter uses governing equations to predict system responses, allowing
control designs to be evaluated without additional physical experiments
\citep{Morris2011,Morris2020}.  We adopt the model-based approach
because, in addition to this practical advantage, the governing equations
provide a basis for analyzing control performance and deriving theoretical
guarantees under explicit assumptions.

Within a model-based approach, actuator design involves two choices: the
spatial configuration of the inputs and the allocation of available
resources among them.  Previous work has treated actuator locations under
linear-quadratic (LQ), $\mathcal H_2$, and $\mathcal H_\infty$ criteria
\citep{ArmaouDemetriou2006,Morris2011,KasinathanMorris2013,MorrisDemetriouYang2015},
as well as actuator supports and shapes
\citep{PrivatTrelatZuazua2017,KaliseKunischSturm2018,EdalatzadehMorris2019}.
Even when the actuator locations and spatial profiles are prescribed,
resource allocation can improve control performance by changing the relative
input strengths \citep{SatoTerasaki2024}.  We study this allocation problem
with the actuator functions fixed and a prescribed total budget.

Controllability allocation scores (CASs) were introduced for
finite-dimensional linear systems as optimal solutions to convex
optimization problems that allocate a fixed input-gain budget among
prescribed input directions for selected target outputs
\citep{sato2026target}.  The volumetric CAS (V-CAS) maximizes the volume of
output changes reachable with unit control energy, whereas the
average-energy CAS (AE-CAS) minimizes the mean energy required for unit
output changes.  CASs generalize the controllability scores developed for
full-state control of finite-dimensional systems with one input per state
variable \citep{SatoTerasaki2024}.  For these earlier scores, existing
theory provides projected-gradient computation, existence, and uniqueness
\citep{SatoTerasaki2024,SatoKawamura2025} and treats time-varying
systems and the infinite-horizon limit
\citep{UmezuSato2025LTV,UmezuSato2026Infinite}.

To formulate the allocation problem for distributed systems governed by
partial differential equations (PDEs), we extend CASs to infinite-dimensional
linear systems with prescribed actuator functions and finitely many target
outputs.  The control objective can specify desired
values for these outputs without prescribing the entire distributed state
\citep{KreindlerSarachik1964,Morse1971,DanhaneLoheacJungers2023}.
For the room example, selecting average temperatures of regions of interest
as outputs focuses the design on those areas.  With $r$ target outputs,
each actuator-to-target response is an $r$-dimensional vector even though
the state remains infinite-dimensional.  Integrating the outer product of
each response gives one $r\times r$ Gramian for each candidate input.
Their budget-weighted sum gives a Gramian for evaluating target controllability
under an allocation and yields finite-dimensional convex allocation problems.

Computing the actuator-to-target responses of these PDE models generally
requires spatial discretization.  The resulting spatially discretized approximations
alter these Gramians and hence the allocation problem, and spatial discretization need
not preserve the observation and control properties of the underlying
infinite-dimensional system
\citep{Zuazua2005}.  Related approximation analyses establish convergence of
optimal costs and actuator locations for LQ design under suitable
assumptions \citep{Morris2011,WuJacobElbern2016}.  We study how
response errors affect target controllability (output controllability)
\citep{BertramSarachik1960} and CASs by tracing them through Gramian assembly
and optimization.

The main contributions are as follows.
\begin{enumerate}
  \item
        We formulate finite-time CASs for infinite-dimensional linear systems with
        bounded distributed input functions and bounded finite-dimensional
        target outputs.  The resulting Gramian depends affinely on the simplex
        allocation, giving convex problems for computing V-CAS and AE-CAS
        after assembling the actuator-to-target responses once.

  \item
        Assuming response convergence, we establish uniform Gramian convergence
        and analyze target controllability under spatial discretization.
        If the continuum system is target controllable, sufficiently fine
        discretizations are also target controllable, and the optimal values
        of the allocation problems converge.  If the allocation problem for
        the continuum system has a unique CAS, the CASs obtained from the
        spatially discretized problems also converge to it.  Under the same
        controllability condition, we obtain an explicit bound on the error
        between these CASs in terms of Gramian
        approximation errors when the Gramian map for the continuum system
        distinguishes allocations.

  \item
        We demonstrate how CASs are used in a two-dimensional
        advection--diffusion example.  Optimal active actuators change with
        circulation speed, while target controllability changes with the horizon
        and can be falsely indicated by coarse upwind discretization.
        Spatial refinement supports the theoretical convergence results.
\end{enumerate}

This paper is organized as follows.
Section~\ref{sec:model-problem} formulates the allocation problem and records
the target-controllability and convexity results used below.
Section~\ref{sec:method} presents the projected-gradient method and numerical
realization.  Section~\ref{sec:theoretical-analysis} establishes spatial
semidiscretization consistency and, when the Gramian map for the continuum
system distinguishes allocations, quantitative allocation-error bounds.
Section~\ref{sec:numerical-experiments}
presents the numerical examples, and Section~\ref{sec:conclusion} concludes
the paper.  Proofs of the results in Section~\ref{sec:theoretical-analysis}
are collected in the \ref{app:theoretical-proofs}.

\paragraph{Notation}

Throughout the paper, $\R^n$ is equipped with the Euclidean inner product and
norm $\|\cdot\|$.
We use $\|\cdot\|_2$ and $\|\cdot\|_{\rm F}$ for the spectral and
Frobenius matrix norms, respectively.
The space of real symmetric $r\times r$ matrices is denoted by $\mathbb S^r$,
and $\mathbb S_{++}^r$ denotes its positive-definite cone.  The symbol $I_n$
denotes the $n\times n$ identity matrix, whereas $I$ denotes the identity
operator on the space understood from context.  We write
$M\succeq0$ and $M\succ0$ for positive semidefiniteness and positive
definiteness.  For $M\in\mathbb S^r$, $\lambda_{\min}(M)$ and
$\lambda_{\max}(M)$ denote its smallest and largest eigenvalues.
For a real Hilbert space $\cH$,
$\langle\cdot,\cdot\rangle_{\cH}$ and $\|\cdot\|_{\cH}$ denote its inner
product and norm.  For matrix-valued functions, $L^2$
norms use the Frobenius norm pointwise.
For a bounded linear operator $L$, $\|L\|$, $L^*$,
$\operatorname{Ran}(L)$, and $\ker(L)$ denote its induced operator norm,
adjoint, range, and kernel, respectively; $\mathcal D(L)$ denotes the domain
when $L$ is unbounded.
For a set $E\subset\R^m$, $\operatorname{ri}(E)$ denotes its relative
interior.  For a nonempty closed convex set $E$, $\Pi_E$ is the Euclidean
projection onto $E$.

\section{Problem Formulation}
\label{sec:model-problem}

We first define the distributed system, its Gramian, and the allocation
problem.
We then characterize target controllability and prove convexity of the allocation
objectives.

\subsection{Controllability allocation problem}
\label{sec:allocation-problem}

Let $\cH$ be a real Hilbert space.  Let $A$ be a linear operator on $\cH$
with domain $\mathcal D(A)\subset\cH$, and assume that $A$ generates a strongly
continuous semigroup $S(t)$ on $\cH$.  Thus,
$S(0)=I$, $S(t+s)=S(t)S(s)$, and $S(t)x\to x$ in $\cH$ as $t\downarrow0$
for every $x\in\cH$.  The semigroup gives the homogeneous evolution
$x(t)=S(t)x(0)$.  Standard realizations of linear diffusion, inflow advection, and
advection--diffusion equations fit this setting
\citep{CurtainZwart1995,TucsnakWeiss2009}.  The numerical example below uses
an advection--diffusion equation and its zero-diffusion limit.
Fix integers $m\geq2$ and $r\geq1$.  We study the state and target-output
equations
\begin{equation}
  \dot x(t)=Ax(t)+B(p)u(t),
  \qquad
  y(t)=Cx(t),
  \label{eq:abstract-system}
\end{equation}
where $x(t)\in\cH$ is the state, $u(t)\in\R^m$ is the control input, and
$y(t)\in\R^r$ is the target output.  Here $A$ is the generator introduced
above.  The vector $p$ is the allocation design variable, and the control $u$
is applied after $p$ is selected.  We call \eqref{eq:abstract-system} the
continuum system to distinguish it from its spatially discretized
approximations.

For a scalar command $u_i(t)$, the fixed function $b_i\in\cH$ produces the
distributed forcing $b_i u_i(t)$ and thereby specifies the spatial action and
nominal gain of candidate actuator $i$.  The allocation-dependent operator
$B(p)$ is built from $b_1,\ldots,b_m$.  We allocate a fixed budget over the standard simplex
\[
  \Delta_m
  :=
  \left\{
  p\in\R^m:
  p_i\geq0\ \text{for all }i,
  \quad
  \sum_{i=1}^m p_i=1
  \right\}.
\]
Its relative interior consists of the allocations with $p_i>0$ for every $i$.
For $p=(p_1,\ldots,p_m)\in\Delta_m$, the component $p_i$ is the fraction of
the squared input-gain budget assigned to channel $i$.  The corresponding
weighted input operator $B(p):\R^m\to\cH$ is
\begin{equation}
  B(p)u
  :=
  \sum_{i=1}^m \sqrt{p_i}\,b_i u_i,
  \qquad
  u=(u_1,\ldots,u_m)\in\R^m.
  \label{eq:weighted-input-operator}
\end{equation}

The $j$th target output is a bounded linear functional of the state.  By the
Riesz representation theorem, it can be written uniquely as
$\langle x,c_j\rangle_{\cH}$ for some $c_j\in\cH$.  The function $c_j$ is the
spatial weighting of target output $j$.  For $c_1,\ldots,c_r$, define the
target operator $C:\cH\to\R^r$ by
\begin{equation}
  Cx
  =
  \begin{pmatrix}
    \langle x,c_1\rangle_{\cH} \\
    \vdots                     \\
    \langle x,c_r\rangle_{\cH}
  \end{pmatrix}.
  \label{eq:target-operator}
\end{equation}
The functions $c_j$ need not coincide with the functions $b_i$, and the
numbers $m$ and $r$ need not agree.  The generator $A$, the two finite
families, and the resulting target operator $C$ are fixed model data; only
$B(p)$ varies with the allocation.

Because the two families are finite, $B(p)$ and $C$ are bounded linear
operators.  Uniformly over $p\in\Delta_m$,
\begin{equation}
  \|B(p)\|
  \leq
  \left(\sum_{i=1}^m p_i\|b_i\|_{\cH}^2\right)^{1/2}
  \leq
  \max_{1\leq i\leq m}\|b_i\|_{\cH},
  \label{eq:input-operator-bound}
\end{equation}
and
\begin{equation}
  \|C\|
  \leq
  \left(\sum_{j=1}^r\|c_j\|_{\cH}^2\right)^{1/2}.
  \label{eq:target-operator-bound}
\end{equation}
The scale of each $b_i$ is part of the actuator model; a common normalization
makes the components of $p$ directly comparable as resource fractions.
When $\cH=\R^n$ and $C$ has full row rank, these definitions give the
finite-dimensional allocation problem of \tcite{sato2026target} with candidate input
directions $b_1,\ldots,b_m$, output matrix $C$, and the full simplex as the
allocation set.

\paragraph{Finite-time Gramian}
Fix a horizon $T>0$.  The mild solution of
\eqref{eq:abstract-system} satisfies
\[
  x(T)
  =
  S(T)x(0)
  +
  \int_0^T S(T-t)B(p)u(t)\dd t.
\]
Hence, the control-induced displacement of the terminal target output is
\[
  y(T)-CS(T)x(0)
  =
  \mathcal R_{p,T}u,
\]
where the target reachability operator $\mathcal R_{p,T}:L^2(0,T;\R^m)\to\R^r$
is defined by
\begin{equation}
  \mathcal R_{p,T}u
  :=
  \int_0^T CS(T-t)B(p)u(t)\dd t.
  \label{eq:target-reachability-operator}
\end{equation}
We call system \eqref{eq:abstract-system} \emph{target controllable at time
  $T$} if there exists $p\in\Delta_m$ such that
\[
  \operatorname{Ran}(\mathcal R_{p,T})=\R^r.
\]
This requirement concerns only the prescribed target output;
it does not require exact controllability of the full state in $\cH$.

The finite-time target controllability Gramian $G(p,T)$ for the continuum
system under allocation $p$ is constructed from the actuator-to-target
responses as follows.
\begin{equation}
  \begin{aligned}
    G(p,T)
     & :=
    \mathcal R_{p,T}\mathcal R_{p,T}^* \\
     & =
    \int_0^T
    CS(t)B(p)B(p)^*S(t)^*C^*
    \dd t
    \in\R^{r\times r}.
  \end{aligned}
  \label{eq:output-gramian}
\end{equation}
Let $B_0:\R^m\to\cH$ be the unweighted actuator operator,
\[
  B_0u:=\sum_{i=1}^m b_i u_i,
\]
and define the response matrix
\begin{equation}
  Q(t)
  :=
  CS(t)B_0
  =
  \begin{bmatrix}
    q_1(t) & \cdots & q_m(t)
  \end{bmatrix}
  \in\R^{r\times m},
  \label{eq:response-matrix}
\end{equation}
where
\[
  q_i(t):=CS(t)b_i.
\]
The Gramian associated with the $i$th candidate input is
\begin{equation}
  G_i(T)
  :=
  \int_0^Tq_i(t)q_i(t)^\top\dd t
  \in\R^{r\times r}.
  \label{eq:channel-gramian-continuous}
\end{equation}
Since
$B(p)=B_0\diag(\sqrt{p_1},\ldots,\sqrt{p_m})$, the Gramian $G(p,T)$ has the
affine representation
\begin{equation}
  G(p,T)=\sum_{i=1}^m p_iG_i(T).
  \label{eq:weighted-gramian}
\end{equation}
Although the input operator $B(p)$ is needed only for $p\in\Delta_m$, we use
\eqref{eq:weighted-gramian} to extend $G(p,T)$ affinely to all $p\in\R^m$
when discussing its positive-definite domain.  Thus, the infinite-dimensional
dynamics enter the allocation problem only through the finite collection of
Gramians $G_1(T),\ldots,G_m(T)$.

As shown in Section~\ref{sec:existing-results}, $G(p,T)\succ0$ is
equivalent to target controllability at time $T$ under the allocation $p$.
We therefore define
\[
  X_T
  :=
  \left\{
  p\in\R^m:
  G(p,T)\succ0
  \right\}.
\]
The set $\Delta_m\cap X_T$ consists of admissible allocations for
which every target direction can be reached with finite input energy.

\paragraph{Controllability allocation scores}
We use the volumetric and average-energy criteria from
controllability scoring \citep{SatoTerasaki2024,sato2026target}.  Here
$G(p,T)$ is assembled from the Hilbert-space actuator-to-target responses in
\eqref{eq:response-matrix}--\eqref{eq:weighted-gramian}.
Subscripts ${\rm V}$ and ${\rm A}$ identify the volumetric and average-energy
criteria, respectively.
For $p\in X_T$, define
\begin{subequations}
  \label{eq:scores}
  \begin{align}
    f_{\rm V,T}(p)&:=-\log\det G(p,T),
    \label{eq:volumetric-objective}\\
    f_{\rm A,T}(p)&:=\tr\bigl(G(p,T)^{-1}\bigr).
    \label{eq:average-energy-objective}
  \end{align}
\end{subequations}
The logarithm avoids direct evaluation of a determinant that may be small,
and the negative sign expresses both criteria as minimization problems.

For either $f_T\in\{f_{\rm V,T},f_{\rm A,T}\}$, the allocation problem for the
continuum system is
\begin{equation}
  \begin{aligned}
     & \text{minimize}   &  & f_T(p)                \\
     & \text{subject to} &  & p\in\Delta_m\cap X_T.
  \end{aligned}
  \label{eq:target-score-problem}
\end{equation}
For $s\in\{{\rm V},{\rm A}\}$, let $p_s^\star(T)$ denote any minimizer of
\eqref{eq:target-score-problem} with $f_T=f_{s,T}$; the argument $(T)$ may be
omitted when the horizon is fixed.  We call such a minimizer a
\emph{controllability allocation score} (CAS) for
problem~\eqref{eq:target-score-problem} under the chosen criterion.
The terms V-CAS and AE-CAS specify the volumetric and average-energy
choices, respectively.

A larger component assigns a larger fraction of the squared input-gain budget
to that candidate, while a zero component means that the candidate is inactive
at the reported optimum.
For general distributed inputs and outputs, CASs allocate resources among
input channels rather than assign centralities to nodes.
If $\Delta_m\cap X_T=\emptyset$, then no admissible allocation makes
system \eqref{eq:abstract-system} target controllable at $T$,
and problem~\eqref{eq:target-score-problem} is infeasible.

\begin{remark}[Geometric and energy interpretation]
  V-CAS maximizes reachable target volume, whereas AE-CAS minimizes the average
  minimum input energy over unit target directions.
  In fact, for $p\in\Delta_m\cap X_T$, the input-energy constraint
  $\int_0^T\|u(t)\|_{\R^m}^2\dd t\leq1$ gives the ellipsoid of reachable
  target displacements
  \[
    \left\{y\in\R^r:y^\top G(p,T)^{-1}y\leq1\right\}.
  \]
  Its volume is proportional to $\sqrt{\det G(p,T)}$, so maximizing it is
  equivalent to minimizing $f_{\rm V,T}(p)$ in
  \eqref{eq:volumetric-objective}.  Moreover, the
  average minimum input energy over unit target directions is proportional to
  $f_{\rm A,T}(p)$ in \eqref{eq:average-energy-objective}.
\end{remark}

\subsection{Target controllability and convexity}
\label{sec:existing-results}

The Hilbert-space target system above differs from the finite-dimensional
full-state model studied by \tcite[Lemmas~3--6 and
  Theorems~4--5]{sato2026projected}.
We first connect target controllability to the Gramian $G(p,T)$ in
\eqref{eq:output-gramian}.  The following proposition characterizes target controllability
under a fixed allocation and the existence of an allocation with this property.
It also characterizes feasibility of
problem~\eqref{eq:target-score-problem}.

\begin{proposition}[Target controllability; adapted from
      {\citep{KreindlerSarachik1964,Morse1971,DanhaneLoheacJungers2023}} and
      {\citep[Lemma~3]{sato2026projected}}]
  \label{prop:target-controllability-feasibility}
  Fix $T>0$.
  \begin{enumerate}
    \item For $p\in\Delta_m$, system \eqref{eq:abstract-system} is target
          controllable at time $T$ under the allocation $p$ if and only if
          $G(p,T)\succ0$.  In that case, the minimum energy required to
          produce $y\in\R^r$ is
          \begin{equation}
            \min_{\mathcal R_{p,T}u=y}
            \|u\|_{L^2(0,T;\R^m)}^2
            =y^\top G(p,T)^{-1}y.
            \label{eq:minimum-target-energy}
          \end{equation}
    \item System \eqref{eq:abstract-system} is target controllable at time
          $T$ if and only if
          \begin{equation}
            \sum_{i=1}^m G_i(T)\succ0.
            \label{eq:aggregate-gramian}
          \end{equation}
          Equivalently, problem~\eqref{eq:target-score-problem} is feasible,
          meaning that
          $\Delta_m\cap X_T\neq\emptyset$.  When this condition holds,
          $\operatorname{ri}(\Delta_m)\subset X_T$; in particular, the system
          is target controllable under the uniform allocation.
  \end{enumerate}
\end{proposition}

\begin{proof}
  The bounded target reachability operator in
  \eqref{eq:target-reachability-operator} satisfies
  $\mathcal R_{p,T}\mathcal R_{p,T}^*=G(p,T)$.
  Since its codomain is $\R^r$, surjectivity is equivalent to
  $G(p,T)\succ0$; the minimum-norm control
  $u^*=\mathcal R_{p,T}^*G(p,T)^{-1}y$ gives
  \eqref{eq:minimum-target-energy}.
  Part~2 follows from the positive-semidefinite matrix argument in
  \tcite[Lemma~3]{sato2026projected}, with allocation set $\Delta_m$
  and Gramians $G_i(T)$.
\end{proof}

\begin{remark}[Horizon dependence]
  \label{rem:horizon-dependence}
  For a finite-dimensional linear time-invariant system and a fixed
  allocation, analyticity of the matrix exponential makes the null space
  $\ker G(p,T)$ independent of $T>0$.  For a strongly continuous semigroup
  with finite propagation, however, $\ker G(p,T)$ can depend on $T$.  For example, on the one-dimensional interval $(0,L)$, let
  $\cH=L^2(0,L)$ and let $S(t)$ be the right-shift semigroup, defined by
  $(S(t)f)(\xi)=f(\xi-t)$ for $\xi>t$ and zero otherwise.  Choose
  $b=\mathbf 1_{(0,a)}$ and $c=\mathbf 1_{(d,d+a)}$, where
  $0<a<d$ and $d+a<L$.  For the scalar input and target, the response
  $q(t)=\langle S(t)b,c\rangle_{\cH}$ vanishes for $0\leq t\leq d-a$ and
  is positive immediately afterward.  Hence the scalar Gramian for this system
  is $G(T)=\int_0^Tq(t)^2\dd t$.  It satisfies $G(T)=0$ for $T\leq d-a$ and
  $G(T)>0$ for $T>d-a$.  Its critical horizon is therefore
  $\inf\{T>0:G(T)>0\}=d-a>0$.  Once target controllability holds,
  $G(p,T_2)\succeq G(p,T_1)$ for $T_2\geq T_1$ ensures that it persists at
  every later horizon.
\end{remark}

\paragraph{Gramian map}
Define the Gramian map for the continuum system by
\begin{equation}
  \Phi_T:\R^m\to\mathbb S^r,
  \qquad
  \Phi_T(v):=\sum_{i=1}^m v_iG_i(T).
  \label{eq:target-gramian-map}
\end{equation}
Then $G(p,T)=\Phi_T(p)$ and
$X_T=\Phi_T^{-1}(\mathbb S_{++}^r)$.  Differences of simplex allocations lie
in
\begin{equation}
  \mathcal L
  :=\{v\in\R^m:\bm1^\top v=0\},
  \label{eq:feasible-direction-space}
\end{equation}
where $\bm1$ is the all-ones vector and $\dim\mathcal L=m-1$.

\begin{remark}[Scope of the matrix analysis]
  Problem~\eqref{eq:target-score-problem} depends on the continuum dynamics
  only through the positive-semidefinite Gramians $G_1(T),\ldots,G_m(T)$.
  Lemmas~3--5 and Theorem~4 of \tcite{sato2026projected} therefore give the
  convexity and existence statements in
  Proposition~\ref{prop:fixed-gramian-consequences}, with allocation set
  $\Delta_m$ and matrix dimension $r$.  These matrix results do not construct
  the responses $q_i$ from the PDE or establish convergence of their spatial
  approximations.  We address these points in
  Sections~\ref{sec:numerical-methods} and~\ref{sec:semidiscrete-consistency}.
\end{remark}

\begin{proposition}[Convexity and existence of CASs; specialization of
      {\citep[Lemmas~4--5 and Theorem~4]{sato2026projected}}]
  \label{prop:fixed-gramian-consequences}
  The following statements hold for
  $f_T\in\{f_{\rm V,T},f_{\rm A,T}\}$.
  \begin{enumerate}
    \item The domain $X_T$ is open and convex, and $f_T$ is smooth and convex
          on $X_T$.  For $p\in X_T$, $v\in\R^m$,
          $G=G(p,T)$, and $\Delta G_v=\Phi_T(v)$,
          \begin{align}
            v^\top\nabla^2f_{\rm V,T}(p)v
             & =\|G^{-1/2}\Delta G_vG^{-1/2}\|_{\rm F}^2,
            \label{eq:vcs-second-directional}                \\
            v^\top\nabla^2f_{\rm A,T}(p)v
             & =2\tr\left((G^{-1}\Delta G_v)^2G^{-1}\right).
            \label{eq:aecs-second-directional}
          \end{align}
    \item If system \eqref{eq:abstract-system} is target controllable at
          time $T$, every nonempty
          sublevel set
          $\{p\in\Delta_m\cap X_T:f_T(p)\leq\ell\}$, $\ell\in\R$, has
          positive distance from
          $\R^m\setminus X_T$.  The extension by $+\infty$ outside
          $\Delta_m\cap X_T$ is proper, lower semicontinuous, and
          level-bounded; consequently,
          problem~\eqref{eq:target-score-problem} admits a CAS.
  \end{enumerate}
\end{proposition}

\begin{proof}
  Apply the matrix arguments of
  \tcite[Lemmas~4--5 and Theorem~4]{sato2026projected}, with allocation
  set $\Delta_m$, Gramians $G_i(T)$, and matrix dimension $r$.
  For part~2, target controllability gives
  condition~\eqref{eq:aggregate-gramian} in
  Proposition~\ref{prop:target-controllability-feasibility}, and hence
  $\Delta_m\cap X_T\neq\emptyset$.
\end{proof}

\section{Method}
\label{sec:method}

This section describes the allocation algorithm and its numerical realization.
First, we state the projected-gradient method for
problem~\eqref{eq:target-score-problem}, whose data are the Gramians $G_i(T)$
in \eqref{eq:channel-gramian-continuous}.
Next, we discretize space and construct the resulting actuator-to-target
responses and Gramians.

\subsection{Optimization on the simplex}
\label{sec:numerical-optimization}

Once the Gramians $G_1(T),\ldots,G_m(T)$ in
\eqref{eq:channel-gramian-continuous} are available, the Gramian $G(p,T)$ in
\eqref{eq:output-gramian} is assembled through
\eqref{eq:weighted-gramian}, and problem~\eqref{eq:target-score-problem}
involves only matrices of size $r\times r$.  The gradient
components are
\begin{align}
  \frac{\partial f_{\rm V,T}}{\partial p_i}(p)
   & =-\tr\bigl(G(p,T)^{-1}G_i(T)\bigr),
  \label{eq:vcs-gradient}                           \\
  \frac{\partial f_{\rm A,T}}{\partial p_i}(p)
   & =-\tr\bigl(G(p,T)^{-1}G_i(T)G(p,T)^{-1}\bigr).
  \label{eq:aecs-gradient}
\end{align}
The gradients and log-determinant are evaluated from a Cholesky factorization
of $G(p,T)$, without forming $G(p,T)^{-1}$ explicitly.

Algorithm~\ref{alg:pgd} is the projected gradient descent (PGD) method of
\tcite[Algorithm~1]{sato2026projected}, written using $G(p,T)$ and $G_i(T)$
and reproduced here to specify the implementation used in the computations.  Since $\Delta_m\cap X_T$ is generally not closed,
the Euclidean projection is taken onto the simplex
$\Delta_m$, not onto $\Delta_m\cap X_T$.  A projected trial point is first
tested for membership in $X_T$ and only then evaluated in the Armijo
condition.  Thus, the objective is never evaluated when $G(p,T)$ is singular.

The projected-gradient mapping is the optimality residual and
stopping criterion of PGD.  For a step size $\alpha>0$, define
\begin{equation}
  \mathcal G_{\alpha}(p)
  :=
  \frac{1}{\alpha}
  \left(
  p-\Pi_{\Delta_m}
  \bigl(p-\alpha\nabla f_T(p)\bigr)
  \right).
  \label{eq:gradient-mapping}
\end{equation}
For $p\in\Delta_m\cap X_T$, the condition $\mathcal G_\alpha(p)=0$ is
equivalent to the first-order optimality condition and hence to global
optimality, since problem~\eqref{eq:target-score-problem} is convex.  Algorithm~\ref{alg:pgd} stops when
the norm of this mapping at the accepted step size is at most
$\varepsilon_{\rm stop}$.

\begin{algorithm}[htbp]
  \caption{Projected-gradient method for CASs}
  \label{alg:pgd}
  \begin{algorithmic}
    \Require Gramians $G_1(T),\ldots,G_m(T)$, objective
    $f_T\in\{f_{\rm V,T},f_{\rm A,T}\}$, initial point
    $p^{(0)}\in\Delta_m\cap X_T$
    \Require Initial trial step $\bar\alpha>0$, backtracking factor
    $\beta\in(0,1)$, Armijo parameter $c\in(0,1)$, and stopping
    tolerance $\varepsilon_{\rm stop}\geq0$
    \For{$k=0,1,2,\ldots$}
    \State Evaluate $f_T(p^{(k)})$ and
    $g^{(k)}=\nabla f_T(p^{(k)})$
    \State $\alpha\gets\bar\alpha$
    \While{true}
    \State $\widehat p\gets
      \Pi_{\Delta_m}(p^{(k)}-\alpha g^{(k)})$
    \If{$\widehat p\in X_T$}
    \If{$f_T(\widehat p)\leq f_T(p^{(k)})
        +c(g^{(k)})^\top(\widehat p-p^{(k)})$}
    \State $\alpha_k\gets\alpha$
    \State $\widetilde p^{(k)}\gets\widehat p$
    \State \textbf{break}
    \EndIf
    \EndIf
    \State $\alpha\gets\beta\alpha$
    \EndWhile
    \If{$\|\widetilde p^{(k)}-p^{(k)}\|/\alpha_k
        \leq\varepsilon_{\rm stop}$}
    \State \Return $p^{(k)}$
    \EndIf
    \State $p^{(k+1)}\gets\widetilde p^{(k)}$
    \EndFor
  \end{algorithmic}
\end{algorithm}

For either objective,
Proposition~\ref{prop:fixed-gramian-consequences} gives convexity, smoothness,
and compact sublevel sets.  Therefore, from every feasible initialization,
backtracking terminates after finitely many reductions and the iterates remain
in the compact initial sublevel set \citep[Theorem~1]{sato2026projected}.
With $\varepsilon_{\rm stop}=0$, the algorithm either terminates at an optimum
or converges to one, without requiring uniqueness
\citep[Corollary~2]{sato2026projected}.

\subsection{Numerical realization}
\label{sec:numerical-methods}

The Gramians $G_i(T)$ in \eqref{eq:channel-gramian-continuous} are built
from the actuator-to-target responses $q_i(t)=CS(t)b_i$, where $S(t)b_i$ is the solution
of the homogeneous PDE with initial state $x(0)=b_i$.  When these responses
cannot be evaluated analytically, we approximate them by spatially discretizing
the PDE.  Integrating the outer products of the resulting responses gives the
Gramians used to optimize the allocation for the spatially discretized system.  For state-space approximations of distributed-parameter
systems and convergence of their finite-time responses, see
\tcite[Sections~4.1 and~7.3]{Morris2020}.

A spatial discretization, indexed by $h$, replaces the PDE state by a vector
$x_h(t)\in\R^{n_h}$ and the generator $A$ by a matrix
$A_h\in\R^{n_h\times n_h}$.  This notation covers finite-difference,
finite-volume, Galerkin, and spectral discretizations.  To compare discrete
states with elements of $\cH$, choose a method-dependent reconstruction
family $\{\psi_{k,h}\}_{k=1}^{n_h}\subset\cH$ and define
$J_h:\R^{n_h}\to\cH$ by
\begin{equation}
  J_hx_h:=\sum_{k=1}^{n_h}(x_h)_k\psi_{k,h}.
  \label{eq:reconstruction-operator}
\end{equation}
The function $J_hx_h$ is the reconstruction of the discrete state in $\cH$.
The symmetric positive-definite mass matrix
$M_h\in\mathbb S_{++}^{n_h}$ represents the discrete state-space inner
product:
\[
  \langle x_h,z_h\rangle_h:=x_h^\top M_hz_h.
\]

Represent the actuator functions $b_i$ in
\eqref{eq:weighted-input-operator} and the target-weighting functions $c_j$ in
\eqref{eq:target-operator} by vectors $b_{i,h},c_{j,h}\in\R^{n_h}$,
respectively.  Collect these vectors in the fixed column matrices:
\[
  \begin{aligned}
    V_{B,h}&:=[b_{1,h}\ \cdots\ b_{m,h}]\in\R^{n_h\times m}, \\
    V_{C,h}&:=[c_{1,h}\ \cdots\ c_{r,h}]\in\R^{n_h\times r}.
  \end{aligned}
\]
The input and target operators are
\[
  B_h(p):=V_{B,h}\diag(\sqrt{p_1},\ldots,\sqrt{p_m}),
  \qquad
  C_h:=V_{C,h}^\top M_h.
\]
Componentwise, the target operator satisfies
\[
  (C_hx_h)_j
  =c_{j,h}^\top M_hx_h
  =\langle x_h,c_{j,h}\rangle_h,
\]
which approximates
$\langle J_hx_h,c_j\rangle_{\cH}=(CJ_hx_h)_j$.

\begin{remark}[Reconstruction and discrete inner product]
  \label{rem:natural-discrete-inner-product}
  The spatial discretization fixes both the reconstruction family
  $\{\psi_{k,h}\}$ in \eqref{eq:reconstruction-operator} and the mass matrix
  $M_h$.  For a two-dimensional cell-centered finite-volume method on a mesh
  $\mathcal T_h$, the entries of $x_h$ are cell values, and
  \[
    J_hx_h=\sum_{K\in\mathcal T_h}(x_h)_K\mathbf1_K,
    \qquad
    M_h=\diag\bigl(|K|:K\in\mathcal T_h\bigr).
  \]
  Thus $x_h^\top M_hz_h=\langle J_hx_h,J_hz_h\rangle_{L^2}$; on a uniform
  rectangular mesh, $M_h=h_xh_yI_{n_h}$.  A two-dimensional finite-difference
  method instead uses nodal values, an interpolating reconstruction, and
  positive quadrature weights, again giving $M_h=h_xh_yI_{n_h}$ on a uniform
  rectangular grid, up to any boundary-weight modifications.  For a
  two-dimensional Galerkin method, the functions $\psi_{k,h}$ form the trial
  basis and
  \[
    (M_h)_{k\ell}
    :=\langle\psi_{k,h},\psi_{\ell,h}\rangle_{\cH}.
  \]
  A spectral method uses the retained modes in the same formula; in particular,
  $M_h=I_{n_h}$ for orthonormal modes.

\end{remark}

The preceding constructions determine the semidiscrete system with allocation
$p$, its actuator-to-target responses, and its Gramians.  Its state and
target-output equations are
\begin{equation}
  \dot x_h(t)=A_hx_h(t)+B_h(p)u(t),
  \qquad
  y_h(t)=C_hx_h(t),
  \label{eq:semidiscrete-system}
\end{equation}
where $y_h(t)\in\R^r$.  The response matrix of
\eqref{eq:semidiscrete-system} is
\begin{equation}
  Q_h(t)
  :=C_he^{tA_h}V_{B,h}
  =\begin{bmatrix}
     q_{1,h}(t) & \cdots & q_{m,h}(t)
   \end{bmatrix}.
  \label{eq:generic-discrete-response}
\end{equation}
Its $i$th column,
$q_{i,h}(t)=C_he^{tA_h}b_{i,h}$, is the actuator-to-target response
corresponding to $q_i(t)$ in $Q(t)$ from \eqref{eq:response-matrix}.
Each response $q_{i,h}$ in \eqref{eq:generic-discrete-response} defines the
Gramian $G_{i,h}(T)$, and allocation $p$ defines $G_h(p,T)$:
\begin{equation}
  \begin{aligned}
    G_{i,h}(T)
     & :=
    \int_0^Tq_{i,h}(t)q_{i,h}(t)^\top\dd t, \\
    G_h(p,T)
     & :=\sum_{i=1}^m p_iG_{i,h}(T).
  \end{aligned}
  \label{eq:semidiscrete-target-gramian}
\end{equation}

Analogously to the continuum system, we define target controllability for the
semidiscrete system~\eqref{eq:semidiscrete-system} by replacing
$CS(T-t)B(p)$ with $C_he^{(T-t)A_h}B_h(p)$ in the reachability operator
\eqref{eq:target-reachability-operator}.  Proposition~\ref{prop:target-controllability-feasibility}
then shows that the system is target controllable under a fixed allocation
$p\in\Delta_m$ if and only if $G_h(p,T)\succ0$.  Such an allocation exists if
and only if $\sum_{i=1}^mG_{i,h}(T)\succ0$, the counterpart of
\eqref{eq:aggregate-gramian}.

Using the Gramians $G_{i,h}(T)$, define the positive-definite domain and
objectives by
\[
  \begin{aligned}
    X_{T,h}
     & :=\{p\in\R^m:G_h(p,T)\succ0\}, \\
    f_{{\rm V},T,h}(p)
     & :=-\log\det G_h(p,T), \\
    f_{{\rm A},T,h}(p)
     & :=\tr\bigl(G_h(p,T)^{-1}\bigr),
     \qquad p\in X_{T,h}.
  \end{aligned}
\]
For $s\in\{{\rm V},{\rm A}\}$, the allocation problem for the semidiscrete
system is
\begin{equation}
  \begin{aligned}
     & \text{minimize}   &  & f_{s,T,h}(p)                \\
     & \text{subject to} &  & p\in\Delta_m\cap X_{T,h}.
  \end{aligned}
  \label{eq:semidiscrete-target-score-problem}
\end{equation}
Let $p_{s,h}^\star(T)$ denote any minimizer of
\eqref{eq:semidiscrete-target-score-problem}.  We call such a minimizer a CAS
for problem~\eqref{eq:semidiscrete-target-score-problem} under the chosen
criterion.
Algorithm~\ref{alg:pgd} applies to
problem~\eqref{eq:semidiscrete-target-score-problem} once the Gramians
$G_{i,h}(T)$ are assembled; subsequent changes in $p$ require only their
weighted sum, so the PDE responses need not be recomputed.

Here and below, $h\to0$ means convergence along the prescribed family of
spatial discretizations.  We assume the following actuator-to-target response
convergence.

\begin{assumption}[Actuator-to-target response consistency]
  \label{ass:channel-response-consistency}
  For every $i=1,\ldots,m$,
  \begin{equation}
    \varepsilon_{i,h}
    :=
    \|q_{i,h}-q_i\|_{L^2(0,T;\R^r)}
    \longrightarrow0
    \qquad\text{as }h\to0.
    \label{eq:channel-response-consistency}
  \end{equation}
\end{assumption}

Under this assumption, Section~\ref{sec:semidiscrete-consistency} analyzes
how spatial approximation affects target controllability, optimal values,
and CASs.

\begin{remark}[PDE-state and output consistency]
  \label{rem:state-output-consistency}
  For each $i$ and $t\in[0,T]$, set
  $x_i(t):=S(t)b_i$, $x_{i,h}(t):=e^{tA_h}b_{i,h}$, and
  $\widehat x_{i,h}(t):=J_hx_{i,h}(t)$.  The triangle inequality gives
  \begin{align*}
    \|q_{i,h}-q_i\|_{L^2(0,T;\R^r)}
    \leq{} &
    \|(C_h-CJ_h)x_{i,h}\|_{L^2(0,T;\R^r)} \\
           & +
    \|C\|\,\|\widehat x_{i,h}-x_i\|_{L^2(0,T;\cH)}.
  \end{align*}
  Thus Assumption~\ref{ass:channel-response-consistency} follows from the two
  limits
  \begin{subequations}
    \label{eq:state-output-sufficient-conditions}
    \begin{align}
      \|\widehat x_{i,h}-x_i\|_{L^2(0,T;\cH)}
       & \xrightarrow[h\to0]{}0,
      \label{eq:reconstructed-trajectory-consistency} \\
      \|(C_h-CJ_h)x_{i,h}\|_{L^2(0,T;\R^r)}
       & \xrightarrow[h\to0]{}0.
      \label{eq:target-operator-consistency}
    \end{align}
  \end{subequations}
  Only \eqref{eq:target-operator-consistency} requires approximation of the
  target functions.  If $M_h$ represents the $\cH$ inner product of
  reconstructed states exactly,
  then
  \[
    \bigl|\bigl((C_h-CJ_h)x_h\bigr)_j\bigr|
    \leq
    \|J_hx_h\|_{\cH}\,\|J_hc_{j,h}-c_j\|_{\cH}.
  \]
  Thus $J_hc_{j,h}\to c_j$ and a uniform
  $L^2(0,T;\cH)$ bound on $J_hx_{i,h}$ imply
  \eqref{eq:target-operator-consistency}.  Quadrature-based inner products
  require the corresponding quadrature consistency.  For a conforming
  Galerkin or spectral method with an exact mass matrix and orthogonally projected
  target functions, $C_h=CJ_h$ identically.
  No rate is assumed; it depends on the PDE, the spatial scheme, the
  reconstruction, and the regularity of $b_i$ and $c_j$.
\end{remark}

\section{Theoretical Analysis}
\label{sec:theoretical-analysis}
\label{sec:semidiscrete-consistency}

We analyze how spatial response errors affect the Gramians and the resulting
allocation problems.  Unlike the finite-dimensional surrogate analysis in
\tcite[Theorems~3--4]{sato2026target}, our analysis starts from $L^2$ errors in
actuator-to-target responses and allows the Gramian $G(p,T)$ in
\eqref{eq:output-gramian} to be singular.  Sections~\ref{sec:gramian-feasibility-consistency}
and~\ref{sec:optimization-consistency} establish preservation of target
controllability and convergence
of optimal values, as well as convergence of CASs when
problem~\eqref{eq:target-score-problem} has a unique CAS.
Section~\ref{sec:unique-allocations} uses restricted injectivity to strengthen
this convergence to a quantitative allocation-error bound.

Throughout this section, the time integrals defining $G_{i,h}(T)$ in
\eqref{eq:semidiscrete-target-gramian} are assumed to be evaluated exactly,
and each computed CAS is assumed to be an exact minimizer.

\subsection{Gramian consistency and target controllability}
\label{sec:gramian-feasibility-consistency}

We derive Gramian error bounds and use them to analyze target controllability
under spatial refinement.
With the responses $q_i,q_{i,h}$ from \eqref{eq:response-matrix} and
\eqref{eq:generic-discrete-response} and their Gramians from
\eqref{eq:channel-gramian-continuous} and
\eqref{eq:semidiscrete-target-gramian}, define the Gramian map for the
semidiscrete system by
\begin{equation}
  \Phi_{T,h}:\R^m\to\mathbb S^r,
  \qquad
  \Phi_{T,h}(v):=\sum_{i=1}^m v_iG_{i,h}(T).
  \label{eq:semidiscrete-gramian-map}
\end{equation}
and, under Assumption~\ref{ass:channel-response-consistency}, set
\begin{equation}
  \begin{aligned}
    \delta_{i,h}
     & :=\bigl(
    \|q_{i,h}\|_{L^2(0,T;\R^r)}
    +\|q_i\|_{L^2(0,T;\R^r)}
    \bigr)\varepsilon_{i,h},                            \\
    \varepsilon_{G,h}
     & :=\max_{1\leq i\leq m}\delta_{i,h},              \\
    \varepsilon_{\Phi,h}
     & :=\left(\sum_{i=1}^m\delta_{i,h}^2\right)^{1/2}.
  \end{aligned}
  \label{eq:semidiscrete-error-quantities}
\end{equation}

\begin{lemma}[Uniform Gramian consistency]
  \label{lem:uniform-gramian-consistency}
  Suppose that Assumption~\ref{ass:channel-response-consistency} holds.
  For every discretization in the prescribed family and every
  $i=1,\ldots,m$,
  \begin{equation}
    \|G_{i,h}(T)-G_i(T)\|_{\rm F}
    \leq\delta_{i,h}.
    \label{eq:channel-gramian-consistency}
  \end{equation}
  Moreover,
  \begin{equation}
    \sup_{p\in\Delta_m}\|G_h(p,T)-G(p,T)\|_{\rm F}
    \leq\varepsilon_{G,h},
    \label{eq:allocated-gramian-uniform-consistency}
  \end{equation}
  and
  \begin{equation}
    \sup_{\|v\|=1}\|\Phi_{T,h}(v)-\Phi_T(v)\|_{\rm F}
    \leq\varepsilon_{\Phi,h}.
    \label{eq:gramian-map-uniform-consistency}
  \end{equation}
  In addition, as
  $h\to0$, the error quantities satisfy
  \[
    \delta_{i,h}\to0\quad(i=1,\ldots,m),
    \qquad
    \varepsilon_{G,h}\to0,
    \qquad
    \varepsilon_{\Phi,h}\to0.
  \]
\end{lemma}

\noindent\emph{Proof.} See
\ref{app:proof-uniform-gramian-consistency}.

\begin{remark}[Numerical quadrature]
Numerical quadrature can be included in these estimates.  If it gives
positive-semidefinite approximations $\widetilde G_{i,h}(T)$ to the Gramians
$G_{i,h}(T)$, with
$\|\widetilde G_{i,h}(T)-G_{i,h}(T)\|_{\rm F}\leq\eta_{i,h}$, then
\eqref{eq:channel-gramian-consistency} and the triangle inequality give
\[
  \|\widetilde G_{i,h}(T)-G_i(T)\|_{\rm F}
  \leq\delta_{i,h}+\eta_{i,h}.
\]
If $\max_i\eta_{i,h}\to0$, replacing $\delta_{i,h}$ by
$\delta_{i,h}+\eta_{i,h}$ in \eqref{eq:semidiscrete-error-quantities}
preserves the convergence of $\varepsilon_{G,h}$ and
$\varepsilon_{\Phi,h}$.  The target-controllability, optimization, and
restricted-injectivity results below then use these enlarged error bounds.  When $G(p,T)$ is singular, the corresponding eigenvalue upper bounds
also acquire the additional terms
$\sum_i p_i\eta_{i,h}$ for an allocation $p$ and $\sum_i\eta_{i,h}$ for the
sum over all candidates.
\end{remark}

Define
\begin{equation}
  \Gamma_T:=\sum_{i=1}^m G_i(T),
  \qquad
  \Gamma_{T,h}:=\sum_{i=1}^m G_{i,h}(T),
  \label{eq:aggregate-gramians}
\end{equation}
for the continuum system~\eqref{eq:abstract-system} and the semidiscrete
system~\eqref{eq:semidiscrete-system}, respectively.

\begin{theorem}[Target controllability under discretization]
  \label{thm:semidiscrete-feasibility-stability}
  Suppose that Assumption~\ref{ass:channel-response-consistency} holds.
  \begin{enumerate}
    \item If the continuum system~\eqref{eq:abstract-system} is target
          controllable at time $T$,
          then there exists
          $\bar h>0$ such that, for every discretization in the prescribed family
          with $0<h<\bar h$,
          $\Gamma_{T,h}\succ0$.  The semidiscrete systems are therefore
          target controllable at time $T$, and
          \[
            \Delta_m\cap X_{T,h}\neq\emptyset,
            \qquad
            \operatorname{ri}(\Delta_m)\subset X_{T,h}.
          \]
    \item If $p\in\Delta_m$ and $G(p,T)$ is singular, then
          \begin{equation}
            0\leq\lambda_{\min}(G_h(p,T))
            \leq\sum_{i=1}^m p_i\varepsilon_{i,h}^2
            \leq\max_{1\leq i\leq m}\varepsilon_{i,h}^2.
            \label{eq:artificial-smallest-eigenvalue-bound}
          \end{equation}
          If $\Gamma_T$ is singular, then similarly
          \begin{equation}
            0\leq\lambda_{\min}(\Gamma_{T,h})
            \leq\sum_{i=1}^m\varepsilon_{i,h}^2.
            \label{eq:aggregate-artificial-eigenvalue-bound}
          \end{equation}
          The rightmost bounds in
          \eqref{eq:artificial-smallest-eigenvalue-bound} and
          \eqref{eq:aggregate-artificial-eigenvalue-bound} both converge to
          zero as $h\to0$.
  \end{enumerate}
\end{theorem}

\noindent\emph{Proof.} See
\ref{app:proof-semidiscrete-feasibility-stability}.

\begin{remark}[Discretization thresholds]
  The threshold $\bar h$ in part~1 of
  Theorem~\ref{thm:semidiscrete-feasibility-stability} depends on
  $\lambda_{\min}(\Gamma_T)$ and the errors $\varepsilon_{i,h}$.  Without a
  quantitative bound for these errors, the theorem does not provide a
  computable mesh size.  Each occurrence of $\bar h$ denotes a local positive
  threshold; statements combining conclusions use the minimum of the
  corresponding thresholds.
\end{remark}

\subsection{Optimal values and CASs}
\label{sec:optimization-consistency}

When the continuum system is target controllable at time $T$,
\tcite[Theorem~7.33]{RockafellarWets1998} gives convergence of optimal values.
If problem~\eqref{eq:target-score-problem} has a unique CAS,
the CASs for problem~\eqref{eq:semidiscrete-target-score-problem} also
converge.

Because the positive-definite domains $X_T$ and $X_{T,h}$ may differ and the
Gramians $G(p,T)$ and $G_h(p,T)$ in \eqref{eq:output-gramian} and
\eqref{eq:semidiscrete-target-gramian}, respectively, can be singular,
ordinary uniform convergence of the objective functions on $\Delta_m$ is not
available.  We therefore encode infeasibility through
extended-value objectives on $\Delta_m$, for which epi-convergence is adapted
to varying effective domains:
\begin{equation}
  \begin{aligned}
    \bar f_{{\rm V},T}(p)
     & :={}
    \begin{cases}
      -\log\det G(p,T), & G(p,T)\succ0,     \\
      +\infty,          & \text{otherwise},
    \end{cases} \\[1ex]
    \bar f_{{\rm A},T}(p)
     & :={}
    \begin{cases}
      \tr(G(p,T)^{-1}), & G(p,T)\succ0,     \\
      +\infty,          & \text{otherwise}.
    \end{cases}
  \end{aligned}
  \label{eq:continuous-extended-score-objectives}
\end{equation}
Define $\bar f_{{\rm V},T,h}$ and $\bar f_{{\rm A},T,h}$ analogously, with
$G_h(p,T)$ in place of $G(p,T)$.  For application of the variational results,
we extend all these functions to $\R^m$ by assigning $+\infty$ outside
$\Delta_m$.

By the sequential characterization in
\tcite[Proposition~7.2]{RockafellarWets1998}, we write
$f_h\xrightarrow[h\to0]{\mathrm{epi}}f$ and say that $f_h$
\emph{epi-converges} to $f$ along the prescribed discretization family if,
for every sequence $h_n\to0$ in that family, the following two
conditions hold: (i) for every $p_n\to p$ with
$p_n,p\in\Delta_m$,
\begin{equation}
  f(p)\leq\liminf_{n\to\infty}f_{h_n}(p_n),
  \label{eq:epi-liminf-definition}
\end{equation}
and (ii) for every $p\in\Delta_m$, there exists $p_n\in\Delta_m$ with
$p_n\to p$ such that
\begin{equation}
  \limsup_{n\to\infty}f_{h_n}(p_n)\leq f(p).
  \label{eq:epi-recovery-definition}
\end{equation}
This formulation allows the available mesh sizes themselves to form a
discrete set.

\begin{theorem}[Convergence of optimal values and CASs]
  \label{thm:semidiscrete-score-consistency}
  Suppose Assumption~\ref{ass:channel-response-consistency} holds.  For
  $s\in\{{\rm V},{\rm A}\}$, the following statements hold.
  \begin{enumerate}
    \item The extended objectives epi-converge:
          \begin{equation}
            \bar f_{s,T,h}\xrightarrow[h\to0]{\mathrm{epi}}\bar f_{s,T}
            \qquad\text{on }\Delta_m.
            \label{eq:semidiscrete-score-epi-convergence}
          \end{equation}
    \item If the continuum system~\eqref{eq:abstract-system} is target
          controllable at time $T$,
          problem~\eqref{eq:target-score-problem} admits a CAS, and there exists
          $\bar h>0$ such that
          problem~\eqref{eq:semidiscrete-target-score-problem} admits a CAS for every
          discretization in the prescribed family with $0<h<\bar h$.
          For every sequence $(h_n)_{n\geq1}$ in that family with $h_n\to0$,
          \begin{equation}
            \min_{p\in\Delta_m}\bar f_{s,T,h_n}(p)
            \xrightarrow[n\to\infty]{}
            \min_{p\in\Delta_m}\bar f_{s,T}(p).
            \label{eq:optimal-value-consistency}
          \end{equation}
          For each $n$, choose a CAS $p_{s,h_n}^\star$ for
          problem~\eqref{eq:semidiscrete-target-score-problem}, so that
          $p_{s,h_n}^\star\in
          \operatorname*{arg\,min}_{p\in\Delta_m}\bar f_{s,T,h_n}(p)$.
          Every cluster point of $\{p_{s,h_n}^\star\}_{n\geq1}$ is a CAS for
          problem~\eqref{eq:target-score-problem}.
    \item \label{item:unique-allocation-consistency}
          Under the assumptions of part~2, if the CAS $p_s^\star$ for
          problem~\eqref{eq:target-score-problem} is unique, then every
          sequence of CASs selected as in part~2 satisfies
          \begin{equation}
            p_{s,h_n}^\star\xrightarrow[n\to\infty]{}p_s^\star.
            \label{eq:optimal-allocation-consistency}
          \end{equation}
  \end{enumerate}
\end{theorem}

\noindent\emph{Proof.} See
\ref{app:proof-semidiscrete-score-consistency}.

\subsection{Error bounds under restricted injectivity}
\label{sec:unique-allocations}

Theorem~\ref{thm:semidiscrete-score-consistency}(\ref{item:unique-allocation-consistency})
gives qualitative CAS convergence when problem~\eqref{eq:target-score-problem}
has a unique CAS, but it does not quantify how spatial approximation
affects the recommended actuator shares.  Such a robustness estimate is important in
implementation because it distinguishes a robust actuator recommendation
from a mesh-dependent numerical artifact.  We therefore impose restricted
injectivity, which quantifies how strongly the Gramian map $\Phi_T$ in
\eqref{eq:target-gramian-map} separates distinct simplex allocations.  Besides yielding strong convexity and uniqueness, this condition
allows Theorem~\ref{thm:quantitative-allocation-consistency} to bound the
allocation error explicitly in terms of Gramian approximation errors.

\begin{assumption}[Restricted injectivity {\citep[Condition~1]{sato2026projected}}]
  \label{ass:restricted-injectivity}
  The map $\Phi_T$ is injective on the space $\mathcal L$ defined in
  \eqref{eq:feasible-direction-space}:
  \begin{equation}
    \ker(\Phi_T)\cap\mathcal L=\{0\}.
    \label{eq:restricted-injectivity}
  \end{equation}
\end{assumption}

The condition excludes nonzero redistributions of the actuator budget that
leave $G(p,T)$ unchanged.  Equivalently,
$G_1(T),\ldots,G_m(T)$ are
affinely independent in $\mathbb S^r$
\citep[Remark~3]{sato2026projected}.

We use this condition to establish strong convexity of the allocation
objectives.  A differentiable $f$ is $\mu$-strongly convex on
$\Delta_m\cap X_T$ if
\[
  \begin{aligned}
    f(q) & \geq f(p)+\nabla f(p)^\top(q-p)
    +\frac{\mu}{2}\|q-p\|^2,                   \\
         & \hspace{7em}p,q\in\Delta_m\cap X_T.
  \end{aligned}
\]
For the volumetric and average-energy objectives, the Hessian formulas
\eqref{eq:vcs-second-directional}--\eqref{eq:aecs-second-directional} give
\begin{equation}
  v^\top\nabla^2f_T(p)v=0
  \quad\Longleftrightarrow\quad
  \Phi_T(v)=0,
  \label{eq:hessian-kernel}
\end{equation}
for $p\in X_T$ and $v\in\R^m$
\citep[Lemma~4]{sato2026projected}.
Thus restricted injectivity excludes zero-curvature directions within
$\mathcal L$.  Since $q-p\in\mathcal L$, it suffices to bound
$v^\top\nabla^2f_T(p)v$ below for $v\in\mathcal L$.

The moduli in the proposition below depend on two constants.  The restricted
separation constant measures the smallest Gramian change per unit
redistribution, and the other constant bounds the largest eigenvalue over all
allocations:
\begin{equation}
  \chi_T:=\min_{\substack{v\in\mathcal L\\\|v\|=1}}
  \|\Phi_T(v)\|_{\rm F},
  \qquad
  \Lambda_T:=\max_{p\in\Delta_m}\lambda_{\max}(G(p,T)).
  \label{eq:fixed-gramian-constants}
\end{equation}
By compactness of the unit sphere in $\mathcal L$,
Assumption~\ref{ass:restricted-injectivity} holds if and only if $\chi_T>0$.
The maximum defining $\Lambda_T$ is finite by continuity and compactness of
$\Delta_m$.

\begin{proposition}[Strong convexity and uniqueness of CASs; specialization of
      {\citep[Lemma~6 and Theorem~5]{sato2026projected}}]
  \label{prop:fixed-gramian-uniqueness}
  If system \eqref{eq:abstract-system} is target controllable at time $T$
  and Assumption~\ref{ass:restricted-injectivity} holds, then the objectives in
  \eqref{eq:volumetric-objective} and \eqref{eq:average-energy-objective}
  are strongly convex on $\Delta_m\cap X_T$, respectively with
  moduli
  \begin{equation}
    \mu_{\rm V}:=\frac{\chi_T^2}{\Lambda_T^2},
    \qquad
    \mu_{\rm A}:=\frac{2\chi_T^2}{\Lambda_T^3}.
    \label{eq:strong-convexity-constants}
  \end{equation}
  Consequently, problem~\eqref{eq:target-score-problem} has a unique CAS for
  each objective.
\end{proposition}

\begin{proof}
  Applying \tcite[Lemma~6]{sato2026projected} to the Gramians $G_i(T)$ on
  $\Delta_m$ gives \eqref{eq:strong-convexity-constants}.
  Existence follows from
  Proposition~\ref{prop:fixed-gramian-consequences}; strong convexity gives
  uniqueness.
\end{proof}

For the semidiscrete system, define
\begin{equation}
  \chi_{T,h}:=\min_{\substack{v\in\mathcal L\\\|v\|=1}}
  \|\Phi_{T,h}(v)\|_{\rm F},
  \quad
  \Lambda_{T,h}:=\max_{p\in\Delta_m}\lambda_{\max}(G_h(p,T)).
  \label{eq:semidiscrete-restricted-separation-constant}
\end{equation}
\begin{proposition}[Restricted injectivity and curvature under discretization]
  \label{prop:semidiscrete-restricted-injectivity}
  Suppose that Assumption~\ref{ass:channel-response-consistency} holds.
  \begin{enumerate}
    \item The quantities $\chi_{T,h}$ and $\Lambda_{T,h}$ satisfy
          \begin{equation}
            |\chi_{T,h}-\chi_T|
            \leq\varepsilon_{\Phi,h},
            \qquad
            |\Lambda_{T,h}-\Lambda_T|
            \leq\varepsilon_{G,h},
            \label{eq:semidiscrete-constant-convergence-bounds}
          \end{equation}
          and hence $\chi_{T,h}\to\chi_T$ and
          $\Lambda_{T,h}\to\Lambda_T$ as $h\to0$.
    \item If Assumption~\ref{ass:restricted-injectivity} also holds, then
          $\chi_T>0$, and there exists $\bar h>0$ such that
          $\Phi_{T,h}$ is injective on $\mathcal L$ for every discretization
          in the prescribed family with $0<h<\bar h$.
    \item If Assumption~\ref{ass:restricted-injectivity} holds and the
          continuum system~\eqref{eq:abstract-system} is target controllable at
          time $T$, then, after decreasing
          $\bar h$ if necessary, for every discretization in the prescribed
          family with $0<h<\bar h$, the volumetric and average-energy objectives
          $f_{{\rm V},T,h}$ and $f_{{\rm A},T,h}$ are strongly convex on
          $\Delta_m\cap X_{T,h}$ with moduli
          \begin{equation}
            \mu_{{\rm V},h}
            :=\frac{\chi_{T,h}^2}{\Lambda_{T,h}^2},
            \qquad
            \mu_{{\rm A},h}
            :=\frac{2\chi_{T,h}^2}{\Lambda_{T,h}^3}.
            \label{eq:semidiscrete-strong-convexity-moduli}
          \end{equation}
          Problem~\eqref{eq:semidiscrete-target-score-problem} therefore has a
          unique CAS for each objective, and the moduli satisfy
          \begin{equation}
            \mu_{{\rm V},h}\xrightarrow[h\to0]{}\mu_{\rm V},
            \qquad
            \mu_{{\rm A},h}\xrightarrow[h\to0]{}\mu_{\rm A}.
            \label{eq:strong-convexity-moduli-convergence}
          \end{equation}
  \end{enumerate}
\end{proposition}

\noindent\emph{Proof.} See
\ref{app:proof-semidiscrete-restricted-injectivity}.

Lemma~\ref{lem:uniform-gramian-consistency} and strong convexity yield the
following CAS error bound.

\begin{theorem}[Error bounds for CASs]
  \label{thm:quantitative-allocation-consistency}
  Suppose that the continuum system~\eqref{eq:abstract-system} is target
  controllable at time $T$ and that
  Assumptions~\ref{ass:channel-response-consistency} and
  \ref{ass:restricted-injectivity} hold.  There exist
  $\underline{\lambda}_T>0$ and $\bar h>0$ such that, for every discretization
  in the prescribed family with $0<h<\bar h$ and each
  $s\in\{\mathrm V,\mathrm A\}$, the following Gramian bounds hold:
  \begin{equation}
    \begin{aligned}
      G(p_s^\star,T)&\succeq\underline{\lambda}_T I_r,
      &G_h(p_s^\star,T)&\succeq\underline{\lambda}_T I_r, \\
      G(p_{s,h}^\star,T)&\succeq\underline{\lambda}_T I_r,
      &G_h(p_{s,h}^\star,T)&\succeq\underline{\lambda}_T I_r.
    \end{aligned}
    \label{eq:optimal-gramian-lower-bounds}
  \end{equation}
  Define
  \begin{equation}
    \begin{aligned}
      \varepsilon_{\nabla,{\rm V},h}
      & :=\frac{\sqrt r}{\underline{\lambda}_T}\varepsilon_{\Phi,h}
          +\frac{\|\Phi_T\|}{\underline{\lambda}_T^2}\varepsilon_{G,h}, \\
      \varepsilon_{\nabla,{\rm A},h}
      & :=\frac{\sqrt r}{\underline{\lambda}_T^2}\varepsilon_{\Phi,h}
          +\frac{2\|\Phi_T\|}{\underline{\lambda}_T^3}\varepsilon_{G,h}.
    \end{aligned}
    \label{eq:allocation-gradient-error-bounds}
  \end{equation}
  Then the CASs $p_s^\star$ and $p_{s,h}^\star$ for
  problems~\eqref{eq:target-score-problem} and
  \eqref{eq:semidiscrete-target-score-problem}, respectively, satisfy
  \begin{equation}
    \|p_{s,h}^\star-p_s^\star\|
    \leq\frac{\varepsilon_{\nabla,s,h}}{\mu_s},
    \qquad s\in\{\mathrm V,\mathrm A\},
    \label{eq:quantitative-allocation-consistency}
  \end{equation}
  where $\mu_{\rm V}$ and $\mu_{\rm A}$ are given by
  \eqref{eq:strong-convexity-constants}.
\end{theorem}

\noindent\emph{Proof.} See
\ref{app:proof-quantitative-allocation-consistency}.

Because $\varepsilon_{G,h},\varepsilon_{\Phi,h}\to0$,
\eqref{eq:quantitative-allocation-consistency} strengthens the qualitative
convergence in
Theorem~\ref{thm:semidiscrete-score-consistency}(\ref{item:unique-allocation-consistency})
under Assumption~\ref{ass:restricted-injectivity}.

\section{Numerical Experiments}
\label{sec:numerical-experiments}

We study a two-dimensional advection--diffusion model with six fixed
actuator functions and three localized target functions.  We first specify
the PDE and its conservative finite-volume discretization.  The numerical
results then present nominal and velocity-dependent allocations,
horizon-dependent target controllability, and spatial refinement.

\subsection{Room model and compactly supported functions}
\label{sec:room-2d}

Let $\Omega=(0,3)\times(0,2)$ and $\cH=L^2(\Omega)$.  Let $x(t,\xi)$ denote
temperature relative to a fixed reference value.  It satisfies
\begin{equation}
  \begin{array}{@{}l@{}}
    \partial_t x
    =\kappa\Delta x-v\cdot\nabla x-\gamma x
    +\displaystyle\sum_{i=1}^{6}\sqrt{p_i}\,b_i u_i(t), \\
    y_j(t)
    =\displaystyle\int_\Omega x(t,\xi)c_j(\xi)\,\dd\xi,
    \qquad j=1,2,3,                                     \\
    (\kappa\nabla x-vx)\cdot n
    =0 \quad\text{on }\partial\Omega,
    \qquad x(0,\cdot)=0.
  \end{array}
  \label{eq:room-2d-pde}
\end{equation}
Here $n$ is the outward unit normal to $\partial\Omega$, $\kappa\geq0$ is
the diffusivity, and $\gamma\geq0$ is the damping rate.
The allocation $p\in\Delta_6$ distributes the fixed actuation budget among
the six independent input channels.  The prescribed velocity $v$ is defined
by the stream function $\psi$:
\begin{equation}
  \begin{aligned}
    \psi(\xi) & =V_0\sin(\pi\xi_1/3)\sin(\pi\xi_2/2),          \\
    v(\xi)    & =(\partial_{\xi_2}\psi,-\partial_{\xi_1}\psi).
  \end{aligned}
  \label{eq:room-2d-streamfunction}
\end{equation}
By construction, $\nabla\cdot v=0$ and $v\cdot n=0$ on the walls.  The
parameter $V_0\geq0$ sets the velocity scale.  All computations use
$\gamma=0.02$; values of $\kappa$ and $V_0$ are stated with the corresponding
results.  For $\kappa>0$, the boundary condition is homogeneous Neumann; for
$\kappa=0$, the walls are characteristic and no additional normal-derivative
condition is imposed.

\begin{figure}[t]
  \centering
  \includegraphics[width=\columnfigurewidth]{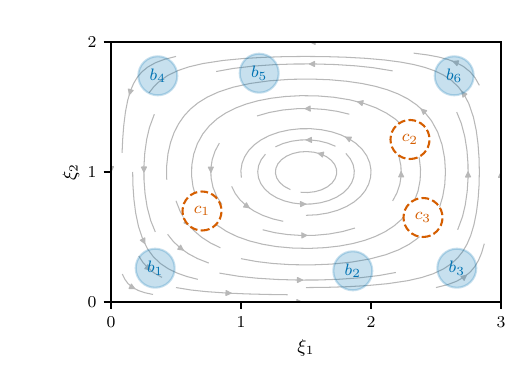}
  \caption{Actuator and target supports in $[0,3]\times[0,2]$.  Filled blue disks show the
    supports of the six actuator functions $b_i$; dashed circles show the
    supports of the three target functions $c_j$.  All radii are $0.15$.
    Gray curves are streamlines of the prescribed velocity field.}
  \label{fig:room-2d-geometry}
\end{figure}

For a center $\zeta$ and radius $\rho=0.15$, set
\begin{equation}
  \beta_\zeta(\xi)
  :=\frac{\sqrt{7/\pi}}{\rho}
  \left(1-\frac{\|\xi-\zeta\|^2}{\rho^2}\right)_+^3,
  \qquad (s)_+:=\max\{s,0\}.
  \label{eq:room-compact-function}
\end{equation}
Each support disk lies inside $\Omega$, and the prefactor gives
$\|\beta_\zeta\|_{L^2(\Omega)}=1$.  Set
\begin{equation}
  \begin{gathered}
    b_i=\beta_{\zeta_i^B},\qquad c_j=\beta_{\zeta_j^C},\\
    \begin{aligned}
      \zeta_1^B & =(0.34,0.26), & \zeta_2^B & =(1.86,0.24), \\
      \zeta_3^B & =(2.66,0.26), & \zeta_4^B & =(0.36,1.74), \\
      \zeta_5^B & =(1.14,1.76), & \zeta_6^B & =(2.64,1.74), \\[2pt]
      \zeta_1^C & =(0.70,0.70), & \zeta_2^C & =(2.30,1.25), \\
      \zeta_3^C & =(2.40,0.65). &           &
    \end{aligned}
  \end{gathered}
  \label{eq:room-input-output-functions}
\end{equation}
All nine supports are mutually disjoint.  In particular, the three target
functions are orthonormal.  Figure~\ref{fig:room-2d-geometry} shows this
configuration and the circulation.  The six actuator supports are not exactly
invariant under the half-turn
$\xi\mapsto(3-\xi_1,2-\xi_2)$.

For fixed $\kappa>0$, the weak Neumann realization of
$A=\kappa\Delta-v\cdot\nabla-\gamma I$ generates a strongly continuous
semigroup on $L^2(\Omega)$.  Since $\nabla\cdot v=0$ and
$v\cdot n=0$, integration by parts gives
\[
  \langle Af,f\rangle_{L^2(\Omega)}
  =-\kappa\|\nabla f\|_{L^2(\Omega)}^2
  -\gamma\|f\|_{L^2(\Omega)}^2.
\]
Thus advection does not change the $L^2$ energy, whereas diffusion and damping
are dissipative.  The smooth first-order term is bounded on $H^1(\Omega)$, so
the standard form construction gives the stated semigroup.

For $\kappa=0$, let $\varphi_t$ denote the characteristic flow of $v$.  This flow
preserves area and $\psi$ and leaves $\Omega$ invariant.  Consequently,
\begin{equation}
  \begin{aligned}
    S(t)f(\xi) & =e^{-\gamma t}f(\varphi_{-t}(\xi)), \\
    (q_i(t))_j & =e^{-\gamma t}\int_\Omega
    c_j(\xi)b_i(\varphi_{-t}(\xi))\,\dd\xi
  \end{aligned}
  \label{eq:room-characteristic-response}
\end{equation}
defines a strongly continuous semigroup and its actuator-to-target
responses.  Hence the homogeneous evolution satisfies the semigroup setting
of Section~\ref{sec:model-problem} for both $\kappa>0$ and $\kappa=0$.

\subsection{Conservative finite-volume discretization}
\label{sec:room-2d-discretization}

Following the standard orthogonal-mesh construction of
\citet[Definition~9.1 and Sections~9.2 and~16]{EymardGallouetHerbin2000}, we
discretize the room by a conservative, cell-centered finite-volume method with
a two-point flux approximation for diffusion and a first-order upwind flux for
advection.  The upwind value is selected according to the sign of the velocity
flux across the face, as in \citet[Section~2]{DelarueLagoutiere2011}.  We use
$N_xN_y$ rectangular cells and set $h_x=3/N_x$, $h_y=2/N_y$, and
$h=\max\{h_x,h_y\}$.  The unknowns are cell values, and we use the
piecewise-constant reconstruction and diagonal mass matrix of cell volumes from
Remark~\ref{rem:natural-discrete-inner-product}.  Since every cell has area
$h_xh_y$, $M_h=h_xh_yI_{n_h}$ and
$x_h^\top M_hx_h=\|J_hx_h\|_{L^2(\Omega)}^2$.  We form the cell averages
of each $b_i$ and $c_j$ by cellwise Gauss quadrature and normalize the
resulting columns in the $M_h$ norm.  These are the columns of $V_{B,h}$
and $V_{C,h}$, respectively; they are not nodal samples.

For an interior face $\sigma$ shared by adjacent cells $K$ and $L$, write
$\sigma=K|L$ and orient its unit normal $n_{K,\sigma}$ outward from $K$.  Let
$F_{K,\sigma}=\int_\sigma v\cdot n_{K,\sigma}\,\dd s$, and let $d_{KL}$ be
the distance between the two cell centers.  We use the total outward
flux
\begin{equation}
  \begin{aligned}
    \mathcal F_{K,\sigma}(x_h)
    ={} & -\kappa\frac{|\sigma|}{d_{KL}}
    \bigl((x_h)_L-(x_h)_K\bigr)                         \\
        & +F_{K,\sigma}^+(x_h)_K+F_{K,\sigma}^-(x_h)_L.
  \end{aligned}
  \label{eq:room-fv-flux}
\end{equation}
Here $F^+=\max(F,0)$ and $F^-=\min(F,0)$.  We evaluate each velocity flux
$F_{K,\sigma}$ exactly from endpoint differences of $\psi$.  These fluxes
therefore sum to zero around every cell, and all wall fluxes vanish.
The generator $A_h$ of the homogeneous semidiscrete
system~\eqref{eq:semidiscrete-system} is defined directly by these fluxes:
\[
  \dot x_h=A_hx_h,
  \qquad
  (A_hx_h)_K
  =-\frac{1}{|K|}\sum_{\sigma\subset\partial K}
  \mathcal F_{K,\sigma}(x_h)-\gamma(x_h)_K.
\]

Interior-face contributions cancel in the total mass balance.  Pairing the
two cell contributions at each interior face shows explicitly how diffusion,
upwinding, and damping dissipate the discrete $L^2$ energy:
\begin{equation}
  \begin{aligned}
    x_h^\top M_hA_hx_h
           ={} & -\kappa\sum_{\sigma=K|L}\frac{|\sigma|}{d_{KL}}
    \bigl((x_h)_L-(x_h)_K\bigr)^2                                \\
                       & -\frac12\sum_{\sigma=K|L}|F_{K,\sigma}|
    \bigl((x_h)_L-(x_h)_K\bigr)^2
    -\gamma x_h^\top M_hx_h.
  \end{aligned}
  \label{eq:room-fv-dissipation}
\end{equation}
Writing $\|x_h\|_h^2:=x_h^\top M_hx_h$, we obtain from
\eqref{eq:room-fv-dissipation}
\begin{equation}
  \|e^{tA_h}x_h\|_h\leq e^{-\gamma t}\|x_h\|_h,
  \qquad t\geq0.
  \label{eq:room-fv-state-bound}
\end{equation}
The second sum in \eqref{eq:room-fv-dissipation} is the numerical diffusion
introduced by first-order upwinding.  It remains on each fixed mesh when
$\kappa=0$ and damps cell-to-cell jumps, producing the coarse-mesh spreading
examined in Section~\ref{sec:room-horizon}.

We compute the responses, CASs, and $\chi_{T,h}$ in
\eqref{eq:semidiscrete-restricted-separation-constant} as follows.  We evaluate
the response matrix in \eqref{eq:generic-discrete-response} using actions of
the matrix exponential.  With only three
targets, propagating $V_{C,h}$ under the adjoint
$A_h^\dagger=M_h^{-1}A_h^\top M_h$ with respect to the $M_h$ inner product
reduces the number of propagated columns.
Composite Simpson quadrature is then applied to the response outer products.
We compute V-CAS and AE-CAS using Algorithm~\ref{alg:pgd}
from the uniform initial allocation, with $\bar\alpha=1$, $\beta=0.5$,
and $c=0.1$.  We set $\varepsilon_{\rm stop}=10^{-12}$ for V-CAS and
$10^{-8}$ for AE-CAS.
To evaluate $\chi_{T,h}$, we use the smallest singular value of the matrix
representing the restriction to $\mathcal L$ of
the Gramian map $\Phi_{T,h}$ in \eqref{eq:semidiscrete-gramian-map}, using
orthonormal coordinates and Frobenius-norm-preserving vectorization of
symmetric matrices.

We next verify the response convergence required by
Assumption~\ref{ass:channel-response-consistency}.  For fixed $\kappa>0$,
\eqref{eq:room-fv-dissipation} bounds the discrete $H^1$ norm in
$L^2(0,T)$, and the finite-volume equations bound the time derivative in
the corresponding dual norm.  Flux consistency and uniqueness of the
PDE solution then give strong convergence of the reconstructed
states by the finite-volume compactness argument
\citep{EymardGallouetHerbin2000}.  For $\kappa=0$, the velocity extends to a
bounded Lipschitz field on $\R^2$, and the zero extensions of $b_i$ are
Lipschitz.  The rectangular refinement meshes have uniformly bounded aspect
ratios, so condition~(4.23) of \citet{DelarueLagoutiere2011} holds.
Their Theorem~5 gives a fully discrete $L^\infty$ error bound that vanishes
as $h\to0$, uniformly over admissible time steps.  Since the wall fluxes
vanish, the whole-space scheme restricts to the room scheme.  On each fixed
mesh, taking the continuous-time limit within the theorem's time-step
restrictions preserves this bound for the evolution operator of the
semidiscrete system.  Damping contributes the factor $e^{-\gamma t}$.  On the
bounded room, these arguments establish
\eqref{eq:reconstructed-trajectory-consistency} for exact actuator cell
averages.  The quadrature-based normalized averages of both actuator and
target functions converge in $L^2(\Omega)$; by
\eqref{eq:room-fv-state-bound}, this preserves the state limit and gives
\eqref{eq:target-operator-consistency}.  Remark~\ref{rem:state-output-consistency}
therefore gives Assumption~\ref{ass:channel-response-consistency} for this
discretization.

\subsection{Allocations with positive diffusion}

First, set $\kappa=0.03$ and $T=6$.  For comparison with the V-CAS and AE-CAS
objectives, define
\[
  W_h(p,T):=\int_0^T e^{tA_h}B_h(p)B_h(p)^\top e^{tA_h^\top}\dd t.
\]
The finite-horizon whole-domain (FH) allocation maximizes
$\tr(M_hW_h(p,T))$ over $p\in\Delta_6$.  Since
\begin{equation}
  \tr\bigl(M_hW_h(p,T)\bigr)
  =\sum_{i=1}^6p_i\int_0^T\|e^{tA_h}b_{i,h}\|_h^2\dd t,
  \label{eq:room-fh-baseline}
\end{equation}
this objective is maximized by setting $p_i=1$ for an actuator with the
largest integral in \eqref{eq:room-fh-baseline} and $p_j=0$ for $j\neq i$.
Spatial
$\mathcal H_2$ measures and modal controllability motivate location design in
\citet{ArmaouDemetriou2006}.  The FH objective in
\eqref{eq:room-fh-baseline} depends on full-state energy, uses no target
information, and restricts the decision to weights on the six
fixed candidates.

Using the Gramian $G_h(p,T)$ in \eqref{eq:semidiscrete-target-gramian},
Table~\ref{tab:room-2d-nominal} compares the selected supports and the values
of $-\log\det G_h$, $\tr G_h^{-1}$, and $\lambda_{\min}(G_h)$ at $V_0=0.15$
on the $191\times127$ mesh.  FH
selects channel 3, V-CAS uses channels $\{1,2,3,6\}$, and AE-CAS uses
$\{1,2,6\}$.  The FH allocation is nearly singular in one target direction and gives larger
values of both $f_{{\rm V},T,h}$ and $f_{{\rm A},T,h}$ than either V-CAS or
AE-CAS.

\begin{table}[t]
  \centering
  \caption{Selected supports and values of $-\log\det G_h$,
    $\tr G_h^{-1}$, and $\lambda_{\min}(G_h)$ at $\kappa=0.03$,
    $V_0=0.15$, and $T=6$ on the $191\times127$ mesh.}
  \label{tab:room-2d-nominal}
  \footnotesize
  \setlength{\tabcolsep}{3pt}
  \begin{tabular}{@{}llrrr@{}}
    \toprule
    Allocation & Support & $-\log\det G_h$
               & $\tr G_h^{-1}$ & $\lambda_{\min}(G_h)$ \\
    \midrule
    uniform & $\{1,\ldots,6\}$ & 26.368 & 29678 & $5.26\times10^{-5}$ \\
    FH      & $\{3\}$            & 40.821 & $6.79\times10^{10}$ & $1.47\times10^{-11}$ \\
    V-CAS   & $\{1,2,3,6\}$      & 25.837 & 30349 & $5.29\times10^{-5}$ \\
    AE-CAS  & $\{1,2,6\}$        & 26.378 & 24251 & $7.99\times10^{-5}$ \\
    \bottomrule
  \end{tabular}
\end{table}

\begin{widefigure}[t]
  \centering
  \includegraphics[width=\linewidth]{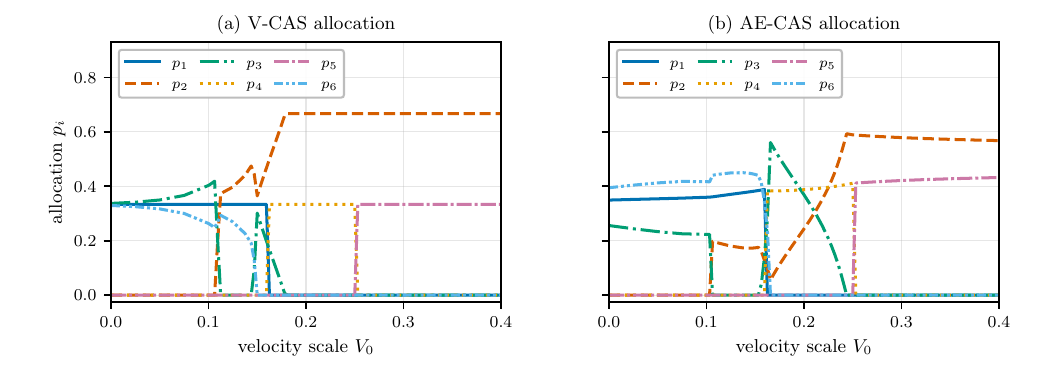}
  \caption{CASs versus $V_0$ at $\kappa=0.03$,
    $\gamma=0.02$, and $T=6$ on the $383\times255$ mesh:
    (a) V-CAS and (b) AE-CAS.  Colors and line styles identify actuator indices;
    lines connect the sampled velocities.}
  \label{fig:room-2d-velocity}
\end{widefigure}

We next examine how the strength of steady advection affects the CASs.  Keeping
$\kappa=0.03$ and $T=6$, we vary $V_0$ over $[0,0.4]$ on the
$383\times255$ mesh.  The PGD solutions have simplex errors at roundoff and
relative first-order residuals below $8\times10^{-12}$.
Figure~\ref{fig:room-2d-velocity} shows the computed CASs.  Stronger
circulation shifts the allocation toward actuators farther
upstream.  At $V_0=0$, V-CAS and AE-CAS both use $\{1,3,6\}$, whose supports
lie near the targets.  The counterclockwise flow makes $b_2$ upstream along the lower
branch and $b_4$ and $b_5$ upstream along the upper-left route.  Accordingly,
V-CAS uses $\{2,4\}$ and AE-CAS uses $\{2,3,4\}$ at $V_0=0.20$;
at $V_0=0.275$, both use the farther-upstream pair $\{2,5\}$, with different
weights.

\subsection{Horizon dependence and apparent target controllability}
\label{sec:room-horizon}

We next set $\kappa=0$ and $V_0=0.15$.  This zero-diffusion transport
setting realizes the finite-propagation mechanism in
Remark~\ref{rem:horizon-dependence} and also illustrates apparent
controllability caused by upwind numerical diffusion.  We classify the
Gramians $G(p,T)$ and $G_h(p,T)$ in \eqref{eq:output-gramian} and
\eqref{eq:semidiscrete-target-gramian}, respectively, as numerically positive
definite when their smallest eigenvalues exceed the practical threshold
$\tau_{\rm PD}:=10^{-14}$.

\begin{widefigure}[t]
  \centering
  \includegraphics[width=\linewidth]{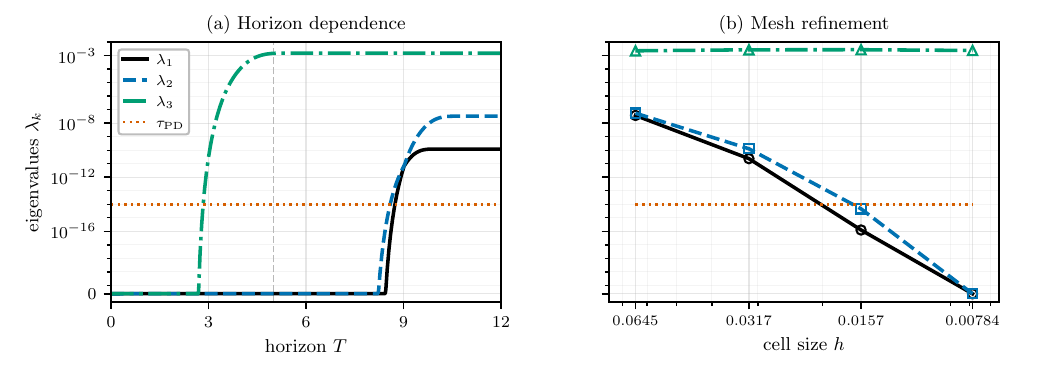}
  \caption{Ordered eigenvalues $\lambda_1\leq\lambda_2\leq\lambda_3$ for
    the uniform allocation, with $\kappa=0$, $V_0=0.15$, and
    $\gamma=0.02$.  (a) Characteristic calculation for the continuum system
    versus horizon $T$; the vertical dashed line marks $T=5$.
    (b) Finite-volume calculation at $T=5$ versus
    $h=\max(3/N_x,2/N_y)$ for meshes $47\times31$, $95\times63$,
    $191\times127$, and $383\times255$.  Colors, line styles, and the
    vertical scale are shared.  Dotted lines show the practical threshold
    $\tau_{\rm PD}$.  The vertical scale is linear near zero and preserves
    the signs of the computed eigenvalues.}
  \label{fig:room-gramian-eigenvalues}
\end{widefigure}

For the continuum system, Figure~\ref{fig:room-gramian-eigenvalues}(a) shows
how the eigenvalues of $G(p,T)$ vary with the horizon under the uniform
allocation.  Since $G(p,T)=\Gamma_T/6$ for this allocation, its positive
definiteness is equivalent to condition~\eqref{eq:aggregate-gramian} in
Proposition~\ref{prop:target-controllability-feasibility}.
Formula~\eqref{eq:room-characteristic-response} provides a reference
independent of the finite-volume operator, which we evaluate using
fourth-order Runge--Kutta integration backward in time and polar quadrature
over the target supports.  At $T=5$, $G(p,5)$ is singular for every allocation.
Invariance of $\psi$ excludes four actuator channels from the first two targets.  For the
remaining two channels, integrating the bound
$|v_2(\xi)|\leq (V_0\pi/3)\sin(\pi\xi_2/2)$ between the support boundaries
gives travel times greater than $5$.  Hence the first two target-response rows
vanish and $\operatorname{rank}G(p,5)\leq1$ for every $p\in\Delta_6$.
At this horizon, the system is not target controllable, and
Assumption~\ref{ass:restricted-injectivity} also fails, since only the $(3,3)$
entries of the six Gramians $G_i(T)$ can be nonzero.
At longer horizons, $G(p,T)$ under the uniform allocation becomes numerically
positive definite, with its smallest eigenvalue crossing the
practical positive-definiteness threshold near $T=8.74$.  This locates the
transition numerically; the bound
$\operatorname{rank}G(p,5)\leq1$ for every $p\in\Delta_6$ does not depend on
that threshold.

Figure~\ref{fig:room-gramian-eigenvalues}(b) shows the finite-volume
eigenvalues at $T=5$ under mesh refinement.  First-order upwinding makes
$G_h(p,5)$ appear positive definite on coarse meshes, but this effect disappears
under refinement.  The values on the two coarsest meshes pass $\tau_{\rm PD}$,
whereas the two finer ones
do not and the finest values are at roundoff level.  The second sum in
\eqref{eq:room-fv-dissipation} spreads responses into targets that are
unreachable in the continuum system over this horizon.  Positive definiteness
on a coarse mesh is therefore not evidence that the continuum system is
target controllable.

\subsection{Spatial refinement and scope of the numerical evidence}
\label{sec:room-2d-validation}

Finally, we examine convergence under spatial refinement and assess numerical
support for target controllability and restricted injectivity.
We fix $\kappa=0.03$, $V_0=0.20$, and $T=6$.  We compare the meshes
$47\times31$, $95\times63$, $191\times127$, and $383\times255$ with a
$767\times511$ reference, using 801 response samples in time.  With
$Q=[q_1\ \cdots\ q_6]$ and the Gramians $G_{i,h}(T)$ in
\eqref{eq:semidiscrete-target-gramian}, define
\begin{equation}
  \begin{aligned}
    E_Q         & :=\frac{\|Q_h-Q_{\rm ref}\|_{L^2(0,T;\R^{3\times6})}}
                    {\|Q_{\rm ref}\|_{L^2(0,T;\R^{3\times6})}},    \\
    E_G         & :=\max_i\frac{\|G_{i,h}-G_{i,\rm ref}\|_{\rm F}}
                          {\|G_{i,\rm ref}\|_{\rm F}},         \\
    E_p^{\rm V} & :=\|p_{{\rm V},h}^\star-p_{{\rm V},\rm ref}^\star\|_\infty, \\
    E_p^{\rm A} & :=\|p_{{\rm A},h}^\star-p_{{\rm A},\rm ref}^\star\|_\infty.
  \end{aligned}
  \label{eq:common-response-error}
\end{equation}

\begin{widefigure}[t]
  \centering
  \includegraphics[width=\linewidth]{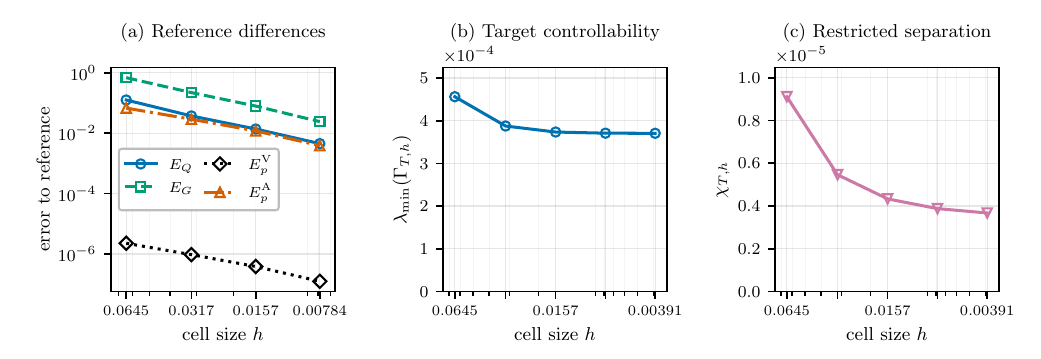}
  \caption{Refinement comparison at $\kappa=0.03$,
    $V_0=0.20$, $\gamma=0.02$, and $T=6$, using 801 time samples.
    (a) $E_Q$, $E_G$, $E_p^{\rm V}$, and $E_p^{\rm A}$
    from \eqref{eq:common-response-error}, relative to the
    $767\times511$ reference mesh; the reference mesh's zero differences
    are omitted.  (b) Smallest eigenvalue of
    $\Gamma_{T,h}=\sum_iG_{i,h}(T)$.
    (c) Computed restricted separation constant $\chi_{T,h}$
    from \eqref{eq:semidiscrete-restricted-separation-constant}.
    Panels (b) and (c) include all five meshes, including the reference,
    and use linear vertical scales starting at zero.
    The horizontal coordinate is $h=\max(3/N_x,2/N_y)$; each panel ends
    at its finest plotted mesh.}
  \label{fig:room-2d-refinement}
\end{widefigure}

The errors $E_Q$, $E_G$, $E_p^{\rm V}$, and $E_p^{\rm A}$ decrease over the
tested meshes in Figure~\ref{fig:room-2d-refinement}(a).  Over these meshes,
$E_Q$ and $E_G$ decay approximately as $h^{1.56}$ and $h^{1.58}$,
respectively, while $E_p^{\rm A}$ decays approximately as $h^{1.35}$.  For V-CAS, every mesh gives
$p_{\rm V}^\star\approx(0,2/3,0,1/3,0,0)$.
Under the assumptions of Theorem~\ref{thm:quantitative-allocation-consistency},
\eqref{eq:allocation-gradient-error-bounds}--\eqref{eq:quantitative-allocation-consistency}
show that if both $\varepsilon_{G,h}$ and $\varepsilon_{\Phi,h}$ are
$O(\varepsilon_h)$ for some $\varepsilon_h\to0$, then
$\|p_{s,h}^\star-p_s^\star\|=O(\varepsilon_h)$: the allocation-error bound
preserves the order of the Gramian approximation errors.  The decay of
$E_G$, $E_p^{\rm V}$, and $E_p^{\rm A}$ in
Figure~\ref{fig:room-2d-refinement}(a) is qualitatively consistent with this estimate.

We use the refinement diagnostics to assess numerical support for target
controllability through the positive-definiteness
condition~\eqref{eq:aggregate-gramian} in
Proposition~\ref{prop:target-controllability-feasibility}, and for restricted
injectivity (Assumption~\ref{ass:restricted-injectivity}).
In Figure~\ref{fig:room-2d-refinement}(b), the smallest eigenvalue of
$\Gamma_{T,h}$ decreases from $4.563\times10^{-4}$ to $3.706\times10^{-4}$
across the five meshes.  The relative change between the two finest meshes
is about $0.08\%$.
Figure~\ref{fig:room-2d-refinement}(c) shows the restricted separation constant
$\chi_{T,h}$ defined in
\eqref{eq:semidiscrete-restricted-separation-constant}.  It measures the
smallest Gramian change per unit allocation change at fixed total weight.
Proposition~\ref{prop:semidiscrete-restricted-injectivity}(1)
gives $\chi_{T,h}\to\chi_T$, with $\chi_T$ from
\eqref{eq:fixed-gramian-constants}; moreover, $\chi_T>0$ is equivalent to
Assumption~\ref{ass:restricted-injectivity}.
The computed values decrease
from $9.110\times10^{-6}$ to $3.672\times10^{-6}$ across the five meshes,
with progressively smaller changes under refinement.
The positivity and stabilization of both diagnostics provide numerical
support for both properties in this parameter setting.

\section{Conclusion}
\label{sec:conclusion}

We formulated finite-time controllability allocation scores for
infinite-dimensional linear systems with bounded distributed actuators and
finite-dimensional target outputs.  The resulting convex allocation problems
use small Gramians $G(p,T)$ assembled from actuator-to-target responses.
Assuming response convergence under spatial refinement, we established uniform
Gramian convergence.
If the continuum system is target controllable at the prescribed final time,
so are the systems on sufficiently fine meshes, and the optimal values of
the allocation problems converge.  CASs also converge when the allocation
problem for the continuum system has a unique CAS.
Under the same controllability condition, restricted injectivity yields strong convexity,
uniqueness, and an explicit allocation-error bound in terms of Gramian
approximation errors.
In the two-dimensional advection--diffusion example, optimal active actuators
changed with circulation speed when diffusion was positive.  Target
controllability depended on the horizon, and coarse upwinding could make
$G_h(p,T)$ appear positive definite even when $G(p,T)$ is singular.  Spatial refinement supported the
theoretical convergence results.

A fully discrete convergence analysis should next incorporate temporal
quadrature errors and optimization residuals alongside spatial response errors.
Validation in real-world settings is another important direction.

\appendix
\begingroup
\renewcommand{\thesection}{Appendix}
\section{Proofs for the Theoretical Analysis}
\label{app:theoretical-proofs}
\endgroup

\subsection{Proof of Lemma~\ref{lem:uniform-gramian-consistency}}
\label{app:proof-uniform-gramian-consistency}

\begin{proof}
  Use the Gramians $G_i(T)$ and $G_{i,h}(T)$ defined in
  \eqref{eq:channel-gramian-continuous} and
  \eqref{eq:semidiscrete-target-gramian}.  Fix $i$ and write
  $\Delta q_{i,h}:=q_{i,h}-q_i$.  Pointwise in time,
  \[
    q_{i,h}q_{i,h}^\top-q_iq_i^\top
    =\Delta q_{i,h}q_{i,h}^\top+q_i\Delta q_{i,h}^\top.
  \]
  Since $\|ab^\top\|_{\rm F}=\|a\|\|b\|$, the triangle inequality and
  the Cauchy--Schwarz inequality give
  \begin{align*}
    \|G_{i,h}(T)-G_i(T)\|_{\rm F}
     & \leq
    \int_0^T\|\Delta q_{i,h}(t)\|
    \bigl(\|q_{i,h}(t)\|+\|q_i(t)\|\bigr)\dd t \\
     & \leq
    \|\Delta q_{i,h}\|_{L^2}
    \bigl(\|q_{i,h}\|_{L^2}+\|q_i\|_{L^2}\bigr)
    =\delta_{i,h}.
  \end{align*}
  For $p\in\Delta_m$,
  \[
    \|G_h(p,T)-G(p,T)\|_{\rm F}
    \leq
    \sum_{i=1}^m p_i\|G_{i,h}(T)-G_i(T)\|_{\rm F}
    \leq\varepsilon_{G,h},
  \]
  where the last inequality uses $p_i\geq0$ and $\sum_i p_i=1$.  Taking the
  supremum proves \eqref{eq:allocated-gramian-uniform-consistency}.  Finally,
  for $\|v\|=1$, the Cauchy--Schwarz inequality gives
  \[
    \|\Phi_{T,h}(v)-\Phi_T(v)\|_{\rm F}
    \leq\sum_{i=1}^m|v_i|\delta_{i,h}
    \leq\left(\sum_{i=1}^m\delta_{i,h}^2\right)^{1/2}
    =\varepsilon_{\Phi,h}.
  \]
  Taking the supremum proves
  \eqref{eq:gramian-map-uniform-consistency}.

  It remains to prove
  $\delta_{i,h}\to0$, $\varepsilon_{G,h}\to0$, and
  $\varepsilon_{\Phi,h}\to0$.  Assumption~\ref{ass:channel-response-consistency}
  gives
  $\|q_{i,h}\|_{L^2}\leq\|q_i\|_{L^2}+\varepsilon_{i,h}$, and hence
  $\delta_{i,h}\to0$ as $h\to0$ for every $i=1,\ldots,m$.  Because the
  number of candidate actuators is finite, the definitions in
  \eqref{eq:semidiscrete-error-quantities} then give
  $\varepsilon_{G,h}\to0$ and $\varepsilon_{\Phi,h}\to0$ as $h\to0$.
\end{proof}

\subsection{Proof of Theorem~\ref{thm:semidiscrete-feasibility-stability}}
\label{app:proof-semidiscrete-feasibility-stability}

\begin{proof}
  Suppose first that the continuum system is target controllable at time $T$.
  Condition~\eqref{eq:aggregate-gramian} in
  Proposition~\ref{prop:target-controllability-feasibility} then gives
  $\lambda_{\min}(\Gamma_T)>0$.  By
  Lemma~\ref{lem:uniform-gramian-consistency},
  \[
    \|\Gamma_{T,h}-\Gamma_T\|_2
    \leq
    \sum_{i=1}^m
    \|G_{i,h}(T)-G_i(T)\|_{\rm F}
    \xrightarrow[h\to0]{}0.
  \]
  Choose $\bar h>0$ so that, for every discretization in the prescribed family
  with $0<h<\bar h$,
  \[
    \|\Gamma_{T,h}-\Gamma_T\|_2
    <\frac12\lambda_{\min}(\Gamma_T).
  \]
  Weyl's inequality \citep[Section~III.2]{Bhatia1997} then yields
  \[
    \lambda_{\min}(\Gamma_{T,h})
    \geq\lambda_{\min}(\Gamma_T)-\|\Gamma_{T,h}-\Gamma_T\|_2
    >\frac12\lambda_{\min}(\Gamma_T)>0.
  \]
  Thus the Gramians $G_{i,h}(T)$ satisfy the counterpart of
  condition~\eqref{eq:aggregate-gramian}.  Part~2 of
  Proposition~\ref{prop:target-controllability-feasibility} gives target
  controllability of the semidiscrete system~\eqref{eq:semidiscrete-system},
  as well as
  $\Delta_m\cap X_{T,h}\neq\emptyset$ and
  $\operatorname{ri}(\Delta_m)\subset X_{T,h}$.

  Fix $p\in\Delta_m$ with $G(p,T)$ singular, and choose a unit
  vector $\eta\in\ker(G(p,T))$.  Since every $G_i(T)$ is positive semidefinite,
  \[
    0=\eta^\top G(p,T)\eta
    =\sum_{i=1}^m p_i\int_0^T|\eta^\top q_i(t)|^2\dd t
  \]
  implies $\eta^\top q_i(t)=0$ for almost every $t$ for every $i$ with
  $p_i>0$.  Therefore
  \begin{align*}
    \eta^\top G_h(p,T)\eta
     & =\sum_{i=1}^m p_i\int_0^T|\eta^\top q_{i,h}(t)|^2\dd t \\
     & =\sum_{i=1}^m p_i
    \int_0^T|\eta^\top(q_{i,h}(t)-q_i(t))|^2\dd t             \\
     & \leq\sum_{i=1}^m p_i\varepsilon_{i,h}^2.
  \end{align*}
  The Rayleigh-quotient characterization and $p\in\Delta_m$ prove
  \eqref{eq:artificial-smallest-eigenvalue-bound}.  If $\Gamma_T$ is
  singular, the argument with a unit vector in $\ker(\Gamma_T)$ gives
  \[
    \lambda_{\min}(\Gamma_{T,h})
    \leq\sum_{i=1}^m\varepsilon_{i,h}^2,
  \]
  proving \eqref{eq:aggregate-artificial-eigenvalue-bound}.
\end{proof}

\subsection{Proof of Theorem~\ref{thm:semidiscrete-score-consistency}}
\label{app:proof-semidiscrete-score-consistency}

\begin{proof}
  We first prove epi-convergence.  Fix an arbitrary sequence $h_n\to0$ in the
  prescribed discretization family, and let $p_n\to p$ in $\Delta_m$.  By
  Lemma~\ref{lem:uniform-gramian-consistency} and continuity of $\Phi_T$,
  \begin{align*}
    \|G_{h_n}(p_n,T)-G(p,T)\|_{\rm F}
     & \leq
    \|G_{h_n}(p_n,T)-G(p_n,T)\|_{\rm F}  \\
     & \quad+\|G(p_n,T)-G(p,T)\|_{\rm F} \\
     & \xrightarrow[n\to\infty]{}0.
  \end{align*}
  If $G(p,T)\succ0$, then $G_{h_n}(p_n,T)\succ0$ for all sufficiently large
  $n$, and continuity of the inverse and log-determinant on
  $\mathbb S_{++}^r$ gives
  $\bar f_{s,T,h_n}(p_n)\to\bar f_{s,T}(p)$.

  It remains to treat singular $G(p,T)$, verify the recovery-sequence condition,
  and extend the result from $\Delta_m$ to $\R^m$.  Along indices for which
  $G_{h_n}(p_n,T)$ is not positive definite, the extended objective is
  $+\infty$.  Along any positive-definite subsequence,
  $\lambda_{\min}(G_{h_n}(p_n,T))\to0$.
  Lemma~\ref{lem:uniform-gramian-consistency} and the definition of $\Lambda_T$
  in \eqref{eq:fixed-gramian-constants} give, for all sufficiently large $n$,
  $\lambda_{\max}(G_{h_n}(p_n,T))\leq\Lambda_T+1$.  Hence
  \[
    \det G_{h_n}(p_n,T)
    \leq(\Lambda_T+1)^{r-1}
    \lambda_{\min}(G_{h_n}(p_n,T))
    \xrightarrow[n\to\infty]{}0,
  \]
  while
  \[
    \tr(G_{h_n}(p_n,T)^{-1})
    \geq\lambda_{\min}(G_{h_n}(p_n,T))^{-1}
    \xrightarrow[n\to\infty]{}+\infty.
  \]
  Thus the epi-liminf inequality
  \eqref{eq:epi-liminf-definition} holds at both positive-definite and singular
  points.  For the epi-limsup inequality, take the constant recovery sequence
  $p_n=p$.  At a positive-definite point,
  Lemma~\ref{lem:uniform-gramian-consistency} gives convergence of the
  objective values; at a singular point, the right-hand side of
  \eqref{eq:epi-recovery-definition} is $+\infty$, so the inequality is
  automatic.  Since the sequence $h_n\to0$ was arbitrary, this proves
  \eqref{eq:semidiscrete-score-epi-convergence} under
  Assumption~\ref{ass:channel-response-consistency} alone.  The extensions also
  epi-converge on $\R^m$: at points outside the closed set $\Delta_m$, every
  convergent sequence eventually remains outside $\Delta_m$, where all
  objectives equal $+\infty$.  At points of $\Delta_m$, terms outside
  $\Delta_m$ cannot lower the liminf, and the same recovery sequence applies.

  For part~2, suppose that the continuum system is target controllable at
  time $T$.  By condition~\eqref{eq:aggregate-gramian} in
  Proposition~\ref{prop:target-controllability-feasibility} and
  Theorem~\ref{thm:semidiscrete-feasibility-stability}, there exists
  $\bar h>0$ such that $\bar f_{s,T}$ and every $\bar f_{s,T,h}$
  with $0<h<\bar h$ are finite at some admissible allocation and never
  take the value $-\infty$, so they are proper.  Continuity on the positive-definite
  domains, divergence at singular Gramians, and closedness of $\Delta_m$
  give lower semicontinuity of their extensions on $\R^m$.  Every sublevel
  set lies in the common compact set $\Delta_m$, so each extension is
  level-bounded.  Thus \tcite[Theorem~1.9]{RockafellarWets1998} gives existence
  of minimizers and finite optimal values, proving the existence assertion
  with this choice of $\bar h$.

  Fix any sequence $h_n\to0$ and discard its finitely many terms, if any,
  for which $\bar f_{s,T,h_n}$ is not proper.  The remaining objectives are
  lower semicontinuous and proper, epi-converge to the proper objective
  $\bar f_{s,T}$, and are eventually level-bounded: for each
  $\ell\in\R$, all sets
  $\{p:\bar f_{s,T,h_n}(p)\leq\ell\}$ lie in $\Delta_m$.
  These are precisely the hypotheses of
  \tcite[Theorem~7.33]{RockafellarWets1998}.  Its conclusions give convergence
  of the optimal values and ensure that every cluster point of selected
  minimizers minimizes $\bar f_{s,T}$ on $\Delta_m$.  Since the sequence
  $h_n\to0$ was arbitrary, part~2 holds for the prescribed family.

  For part~3, uniqueness of the minimizer $p_s^\star$ of $\bar f_{s,T}$ makes every
  cluster point in part~2 equal to $p_s^\star$.  Compactness of
  $\Delta_m$ gives convergence of the whole sequence.
\end{proof}

\subsection{Proof of
  Proposition~\ref{prop:semidiscrete-restricted-injectivity}}
\label{app:proof-semidiscrete-restricted-injectivity}

\begin{proof}
  For every $v\in\mathcal L$ with $\|v\|=1$,
  Lemma~\ref{lem:uniform-gramian-consistency} and the definitions of
  $\chi_T$ and $\chi_{T,h}$ in \eqref{eq:fixed-gramian-constants} and
  \eqref{eq:semidiscrete-restricted-separation-constant} give
  \[
    \bigl|
    \|\Phi_{T,h}(v)\|_{\rm F}-\|\Phi_T(v)\|_{\rm F}
    \bigr|
    \leq\varepsilon_{\Phi,h}.
  \]
  Taking minima over the unit sphere first gives
  $\chi_{T,h}\geq
    \chi_T-\varepsilon_{\Phi,h}$.  Interchanging
  $\Phi_{T,h}$ and $\Phi_T$ gives the reverse inequality, and hence
  $|\chi_{T,h}-\chi_T|
    \leq\varepsilon_{\Phi,h}$.  For every
  $p\in\Delta_m$, Weyl's inequality and
  \eqref{eq:allocated-gramian-uniform-consistency} give
  \[
    \begin{aligned}
       & |\lambda_{\max}(G_h(p,T))-\lambda_{\max}(G(p,T))| \\
       & \qquad\leq\|G_h(p,T)-G(p,T)\|_2
      \leq\varepsilon_{G,h}.
    \end{aligned}
  \]
  Taking maxima over $\Delta_m$ proves the second bound in
  \eqref{eq:semidiscrete-constant-convergence-bounds} and completes part 1.

  Under Assumption~\ref{ass:restricted-injectivity}, compactness of the unit
  sphere in $\mathcal L$ gives $\chi_T>0$.  Since
  $\varepsilon_{\Phi,h}\to0$, choose $\bar h>0$ so that
  $\chi_{T,h}\geq\chi_T/2>0$ whenever $0<h<\bar h$.  This proves part~2.  If
  the continuum system is also target controllable at time $T$, decrease
  $\bar h$ if necessary so that
  Theorem~\ref{thm:semidiscrete-feasibility-stability} gives target
  controllability of the semidiscrete systems whenever $0<h<\bar h$.

  For each such $h$, the semidiscrete system therefore satisfies the
  target-controllability and restricted-injectivity conditions of
  Proposition~\ref{prop:fixed-gramian-uniqueness}.  Applying its curvature
  bounds \eqref{eq:strong-convexity-constants} gives precisely the moduli in
  \eqref{eq:semidiscrete-strong-convexity-moduli}.
  Their convergence follows from part~1, since
  condition~\eqref{eq:aggregate-gramian} implies $\Lambda_T>0$ by
  \eqref{eq:fixed-gramian-constants}.
  Part~2 of Proposition~\ref{prop:fixed-gramian-consequences} applies to these
  Gramians $G_{i,h}(T)$ as well.  Strong convexity gives uniqueness, proving
  part~3.
\end{proof}

\subsection{Proof of Theorem~\ref{thm:quantitative-allocation-consistency}}
\label{app:proof-quantitative-allocation-consistency}

\begin{proof}
  We first bound the smallest eigenvalues of the four Gramians in
  \eqref{eq:optimal-gramian-lower-bounds} uniformly away from zero.  Set
  $p^0=m^{-1}\bm1$.  By condition~\eqref{eq:aggregate-gramian} in
  Proposition~\ref{prop:target-controllability-feasibility}, target
  controllability gives $G(p^0,T)=m^{-1}\Gamma_T\succ0$, hence $p^0\in X_T$.
  By Lemma~\ref{lem:uniform-gramian-consistency}, this point
  also belongs to $X_{T,h}$ for all sufficiently small $h$, and
  $f_{s,T,h}(p^0)\to f_{s,T}(p^0)$ for
  $s\in\{\mathrm V,\mathrm A\}$.  Define
  $C_s:=1+\max\{0,f_{s,T}(p^0)\}$.  These constants are positive and
  independent of $h$, and
  $f_{s,T,h}(p^0)\leq C_s$ for both criteria and all sufficiently small $h$.
  By \eqref{eq:semidiscrete-constant-convergence-bounds}, with $\Lambda_T$ and
  $\Lambda_{T,h}$ from \eqref{eq:fixed-gramian-constants} and
  \eqref{eq:semidiscrete-restricted-separation-constant}, we may also require
  $\Lambda_{T,h}\leq\Lambda_T+1$ by taking $h$ sufficiently small.

  Because $p_{s,h}^\star$ minimizes $f_{s,T,h}$ and $p^0$ is
  feasible, we have
  \begin{equation}
    f_{s,T,h}(p_{s,h}^\star)
    \leq f_{s,T,h}(p^0)\leq C_s.
    \label{eq:proof-optimal-objective-bound}
  \end{equation}
  For the volumetric objective,
  \eqref{eq:proof-optimal-objective-bound} implies that the determinant of
  $G_h(p_{{\rm V},h}^\star,T)$ is at least $e^{-C_{\rm V}}$.  Since each of
  its other $r-1$ eigenvalues is at most $\Lambda_T+1$, we obtain the first
  estimate below.  For the average-energy objective,
  \eqref{eq:proof-optimal-objective-bound} and the inequality between the
  trace and the largest eigenvalue of the inverse give the second.  Thus
  \begin{equation}
    \begin{aligned}
      \lambda_{\min}(G_h(p_{{\rm V},h}^\star,T))
      &\geq \frac{e^{-C_{\rm V}}}{(\Lambda_T+1)^{r-1}}, \\
      \lambda_{\min}(G_h(p_{{\rm A},h}^\star,T))
      &\geq \frac{1}{C_{\rm A}}.
    \end{aligned}
    \label{eq:proof-discrete-optimal-gramian-lower-bounds}
  \end{equation}
  The bounds in \eqref{eq:proof-discrete-optimal-gramian-lower-bounds}
  also hold for $G(p_s^\star,T)$, since
  $f_{s,T}(p_s^\star)\leq f_{s,T}(p^0)\leq C_s$ and
  $\lambda_{\max}(G(p_s^\star,T))\leq\Lambda_T+1$.  Set
  \[
    \underline{\lambda}_T
    :=\frac12\min\left\{
    \frac{e^{-C_{\rm V}}}{(\Lambda_T+1)^{r-1}},
    \frac{1}{C_{\rm A}}
    \right\}.
  \]
  Both $G(p_s^\star,T)$ and $G_h(p_{s,h}^\star,T)$ are then bounded below
  by $2\underline{\lambda}_T I_r$ for each criterion.  Choose $\bar h>0$
  so that the preceding estimates hold and
  $\varepsilon_{G,h}\leq\underline{\lambda}_T$ for every discretization in
  the prescribed family with $0<h<\bar h$.  Applying
  \eqref{eq:allocated-gramian-uniform-consistency} at both optimizers gives
  \eqref{eq:optimal-gramian-lower-bounds}, with this $h$-independent choice
  of $\underline{\lambda}_T$ common to both criteria.

  We next estimate the difference between the gradients of $f_{s,T}$ and
  $f_{s,T,h}$ at $p_{s,h}^\star$.  For $M\succ0$, the matrix
  gradients of $-\log\det M$ and $\tr(M^{-1})$ are $-M^{-1}$ and $-M^{-2}$,
  respectively.  Taking adjoints with respect to the Euclidean and
  Frobenius inner products and applying the chain rule, we obtain
  \begin{equation}
    \begin{aligned}
      \nabla f_{{\rm V},T}(p)
      &=-\Phi_T^*\bigl(G(p,T)^{-1}\bigr), \\
      \nabla f_{{\rm A},T}(p)
      &=-\Phi_T^*\bigl(G(p,T)^{-2}\bigr), \\
      \nabla f_{{\rm V},T,h}(p)
      &=-\Phi_{T,h}^*\bigl(G_h(p,T)^{-1}\bigr), \\
      \nabla f_{{\rm A},T,h}(p)
      &=-\Phi_{T,h}^*\bigl(G_h(p,T)^{-2}\bigr).
    \end{aligned}
    \label{eq:proof-objective-gradient-factorizations}
  \end{equation}
  Whenever $M,N\succeq\underline{\lambda}_T I_r$, the inverse identities imply
  \begin{equation}
    \begin{aligned}
      \|M^{-1}-N^{-1}\|_{\rm F}
      &\leq \underline{\lambda}_T^{-2}\|M-N\|_{\rm F}, \\
      \|M^{-2}-N^{-2}\|_{\rm F}
      &\leq2\underline{\lambda}_T^{-3}\|M-N\|_{\rm F}, \\
      \|M^{-1}\|_{\rm F}
      &\leq\frac{\sqrt r}{\underline{\lambda}_T},
      \qquad
      \|M^{-2}\|_{\rm F}
      \leq\frac{\sqrt r}{\underline{\lambda}_T^2}.
    \end{aligned}
    \label{eq:proof-inverse-perturbation-bounds}
  \end{equation}
  Combining \eqref{eq:proof-objective-gradient-factorizations} and
  \eqref{eq:proof-inverse-perturbation-bounds} with
  \eqref{eq:allocated-gramian-uniform-consistency}--\eqref{eq:gramian-map-uniform-consistency},
  and using the equality of the operator norms of a map and its adjoint, gives
  \begin{equation}
    \|\nabla f_{s,T,h}(p_{s,h}^\star)
       -\nabla f_{s,T}(p_{s,h}^\star)\|
    \leq\varepsilon_{\nabla,s,h},
    \label{eq:proof-optimal-gradient-discrepancy}
  \end{equation}
  where $\varepsilon_{\nabla,s,h}$ is defined in
  \eqref{eq:allocation-gradient-error-bounds}.

  We now use \eqref{eq:proof-optimal-gradient-discrepancy} to bound the
  allocation error.  Fix $s\in\{\mathrm V,\mathrm A\}$ and define
  $d_h:=p_{s,h}^\star-p_s^\star\in\mathcal L$.  By
  \eqref{eq:optimal-gramian-lower-bounds}, both optimizers belong to
  $\Delta_m\cap X_T\cap X_{T,h}$.
  Either optimizer is therefore an admissible comparison point in the other
  problem.  The first-order optimality
  conditions are therefore
  \begin{equation}
    \nabla f_{s,T}(p_s^\star)^\top d_h\geq0,
    \qquad
    \nabla f_{s,T,h}(p_{s,h}^\star)^\top d_h\leq0.
    \label{eq:proof-cross-optimality-conditions}
  \end{equation}
  Strong convexity of $f_{s,T}$ with modulus $\mu_s$, followed by
  \eqref{eq:proof-cross-optimality-conditions},
  \eqref{eq:proof-optimal-gradient-discrepancy}, and the Cauchy--Schwarz
  inequality, gives
  \begin{equation}
    \begin{aligned}
      \mu_s\|d_h\|^2
      &\leq
      \bigl(\nabla f_{s,T}(p_{s,h}^\star)
             -\nabla f_{s,T}(p_s^\star)\bigr)^\top d_h \\
      &\leq\nabla f_{s,T}(p_{s,h}^\star)^\top d_h \\
      &\leq
      \bigl(\nabla f_{s,T}(p_{s,h}^\star)
             -\nabla f_{s,T,h}(p_{s,h}^\star)\bigr)^\top d_h \\
      &\leq\varepsilon_{\nabla,s,h}\|d_h\|.
    \end{aligned}
    \label{eq:proof-allocation-error-chain}
  \end{equation}
  If $d_h\neq0$, dividing \eqref{eq:proof-allocation-error-chain} by
  $\mu_s\|d_h\|$ proves
  \eqref{eq:quantitative-allocation-consistency}.  The same conclusion is
  immediate when $d_h=0$.
\end{proof}

\section*{Acknowledgements}
The authors thank Yuito Nakabe for helpful discussions on controllability
scores for infinite-dimensional systems.

This work was supported by JST PRESTO, Japan, Grant Number \mbox{JPMJPR25K4}.

% Alphabetical author--year references shared by both active drivers.
% The bibliographic records mirror references_shared.tex; only order and style differ.

\end{document}